\documentstyle[amssymb,amscd,10pt]{amsart}
\input xy
\xyoption{all}

\newtheorem{thm}{Theorem}[section]
\newtheorem{cor}[thm]{Corollary}

\newtheorem{prop}[thm]{Proposition}
\newtheorem{lem}[thm]{Lemma}

\theoremstyle{definition}
\newtheorem{dfn}[thm]{Definition}

\newtheorem{rem}[thm]{Remark}
\numberwithin{equation}{section}

\begin{document}

\newcommand{\id}{\mathop{}\mathopen{}\mathrm{id}}
\newcommand{\Hom}{\operatorname{Hom}}
\newcommand{\Spec}{\operatorname{Spec}}
\newcommand{\Aut}{\operatorname{Aut}}
\newcommand{\Ker}{\operatorname{Ker}}
\newcommand{\Img}{\operatorname{Im}}
\newsymbol\cuadrado1003

\title[Automorphism group of a toric variety ]
{Automorphism group of a toric variety}
\author{M.T. Sancho de Salas \\ J.P. Moreno Gonz\'{a}lez\\ C. Sancho de Salas
}
\address{ \newline M.Teresa Sancho de Salas\newline Departamento de
Matem\'{a}ticas\newline
Universidad de Salamanca\newline  Plaza de la Merced 1-4\\
37008 Salamanca\newline  Spain}
\email{sancho@@usal.es}
\address{\newline Jesus Moreno Gonzalez\newline Departamento de
Matem\'{a}ticas\newline
Universidad de Salamanca\newline  Plaza de la Merced 1-4\\
37008 Salamanca\newline  Spain}
\email{jmorenogon@@yahoo.es}
\address{ \newline Carlos Sancho de Salas\newline Departamento de
Matem\'{a}ticas\newline
Universidad de Salamanca\newline  Plaza de la Merced 1-4\\
37008 Salamanca\newline  Spain}
\email{mplu@@usal.es}

\thanks{ The third author was supported by research project MTM2013-45935-P (MINECO).}

\begin{abstract}
We calculate the automorphism group of a complete toric variety $X$ with torus $T_M$. We prove that the radical unipotent of $Aut_k^0X$ is a semidirect product of additive groups, the reductive part is a quotient of a product of lineal groups and we give the action of linear groups on the additive groups. We also prove that $Aut_kX/Aut_k^0X$ is a quotient of the group of automorphisms of $M$ leaving invariant the fan.
\end{abstract}
\setcounter{section}{-1} \maketitle
\tableofcontents

\section{Introduction}

The aim of this article is determined thoroughly the group of automorphisms of a  complete toric variety $X$ with torus $T_M$. We will always work over an algebraically closed field and of characteristic $0$.

The first theorem is to identify $Aut^0_kX$  with a quotient of the group of graded automorphisms of the algebra known as Cox ring (Theorem 5.6). This is proven in [C]  for the case  simplicial.

 It is known that every affine algebraic group is a semidirect product of an unipotent group and a reductive group (Prop 5.4.1 in [Co]). We determine and construct the reductive subgroup and the unipotent subgroup  of $Aut^0_kX$ \ as well as the action of the reductive group on unipotent group.

This description only depends on the rays of the fan of the toric variety. The set of the rays, $\Delta_1$, is identified with the set of irreducible hypersurfaces of $X$ invariants by $T_M$ and we consider in it the equivalence relation given by linear equivalence of divisors. We denote by $W$  the quotient set of $\Delta_1$  modulo this relation.

 For each class of equivalence $F\in W$, we construct an additive subgroup $V_F$ of $Aut^0_kX$ and a linear group $Gl_F$ acting on $X$.  The unipotent radical of $Aut_k^0X$ is obtained as a semidirect product of the groups $V_F$ ordered  thanks to a partial order  defined on $W$ (Theorem 6.14). These additive groups are obtained as one-parameter subgroups associated to the no semisimple roots of the group [D].

 The reductive group is obtained as a quotient of the product of the linear groups $Gl_F$. (theorem 6.6 and 7.2).

 The theorem 6.17 determines the action of the linear groups on the additive groups giving the irreducible representations.

In the last section we describe the finite group \  $Aut_kX/Aut^0_kX$ \  as a quotient of the  group  of the automorphisms of $M$   leaving invariant the fan (theorem 8.8). This theorem is also proved by Cox for the case simplicial.

\section{Group functors}

We will  work with several groups  representing   functors of groups.

a) The automorphism group of a complete variety $X$ over $k$. We denote it by $Aut_kX$. It is the representant of the functor over the category of $k$-schemes given by  \  $Z \ \rightsquigarrow \  Aut_{Z}(X\times Z)$.
 This functor is representable. See  [MO].

b) Given a ${\mathbb Z}$-module $N$ and $k[N]=\{\underset{\alpha\in N}{\sum} \lambda_\alpha x^\alpha \  / \  \lambda_\alpha\in k\text{ and almost all zero }\}$, the multiplicative group \  $T_N=Spec\,k[N]$ \   is the representant of the functor over category of $k$-schemes given by  \  \  $Z \ \rightsquigarrow \  Hom_{gr}(N,{{\mathcal O}}_{Z}(Z)^\times)\text{ (group morphisms)}$
 \  \  where ${{\mathcal O}}_{Z}(Z)^\times$ \  is the group of invertible elements of ${\mathcal O}_{Z}(Z)$.

We will put $\lambda^\beta=\lambda(\beta)$ for each $\lambda\in Hom_{gr}(N,{\mathcal O}_{Z}(Z)^\times)$ and $\beta \in N$.

Thus the  torus \ $T_{{\mathbb Z}^r}$ \  is the representant of the functor:   $Z \  \rightsquigarrow \  {\mathcal O}_{Z}(Z)^{\times  r}$ \  and if  $\lambda=(\lambda_1,\dots,\lambda_r)\in {\mathcal O}_{Z}(Z)^\times{^r}$ \  and \  $m=(n_1,\dots,n_r)\in {\mathbb Z}^r$, then $\lambda^m=\lambda_1^{n_1}\cdot{\dots}\cdot\lambda_r^{n_r}$.

c) Group of automorphisms of a $N$-graded algebra.

Let $A=\underset{s\in N}{\oplus} A_s$ be a $N$-graded algebra.

The group of $ N$-graded automorphisms of $A$ is the scheme  representing the functor over the category of $k$-schemes given by \  $Y \rightsquigarrow Aut_{{\mathcal O}_Y-gra}(A\otimes_k{\mathcal O}_Y)$ (graded automorphisms).

If $A_0=k$ and $A$ is of finite type over $k$, this functor is representable and its representant  is denoted by \  $Aut_{k-gr}A$. This is because if $A$ is generated by $\xi_1,\dots,\xi_r$ and the degree of $\xi_i=s_i$, then $Aut_{k-gr}A$ is a subgroup of $Aut_kV$ being $V=\overset{r}{\underset{i=1}{\bigoplus}} A_{s_i}$ a $k$-vectorial space of finite dimension.

 Because of being $A$ a $ N$-graded algebra,   $T_{ N}=Spec\,k[ N]$ \  is a subgroup of $Aut_{k-gr}A $:
 For each \  $\lambda\in Hom_{gr}(N,{\mathcal O}_{Z}(Z)^\times)$, $\lambda$ produces the \  following ${\mathcal O}_{Z}$-graded automorphism: $h_\lambda: A\otimes_k{\mathcal O}_{Z}\to A\otimes_k{\mathcal O}_{Z}$ \  given by \   $h_\lambda(a_s)=\lambda^s\cdot a_s$ \ for each \  $a_s\in A_s$.

d) The linear group. Let  $V_l$ be a $k$-vectorial space of finite dimension $l$. We denote $Aut_kV_l$ the representant of the group functor such that for each $k$-scheme $Y$ corresponds the group of the automorphisms of  \  $V_l\otimes_k{\mathcal O}_Y$ \  as ${\mathcal O}_Y$-module.

If we choose a basis $(z_1,\dots,z_l)$ on $V_l$, the linear group is the group $Aut_{k-gr}A$ where $A=k[z_1,\dots,z_l]$ is the  ${\mathbb Z}$-graded algebra such that the degree of each $z_i$ is 1.

If $C=(c_{ij})$ is an invertible matrix of order $l$, it defines the ${\mathbb Z}$-graded automorphism \  $\tau_C:k[z_1,\dots,z_l]\to k[z_1,\dots,z_l]$ \  given by $\tau_C(z_j)=\overset{l}{\underset{i=1}{\sum}} c_{ij}z_i$.

e) The quotient group \  $Aut_{k-gr}A\,/\,T_{ N}$ \ where $A$ is a $N$-graded algebra with the above conditions.

$T_{ N}$ is within the center of $Aut_{k-gr}A$ and so \ $Aut_{k-gr}A\,/\,T_{ N}$ is a group [SGA3].

To give a morphism \  $f:Y\to Aut_{k-gr}A\,/\,T_{N}$ \ is sufficient to give an open covering $\{U_i\}_{i\in I}$ on $Y$ and  \ $\tau_i\in Aut_{{\mathcal O}_{U_i}-gr}(A\otimes_k{\mathcal O}_{U_i})$  \  such that for all $i,j\in I$,
 ${\tau_i}_{|U_i\cap U_j}=h_\lambda\circ {\tau_j}_{|U_i\cap U_j}$ \  for some $\lambda\in Hom_{gr}(N,{\mathcal O}_Y(U_i\cap U_j)^\times)$.

f) The additive group of dimensi\'{o}n $r$. It is the underlying  group under addition of   a $k$-vectorial space $V_r$. If we choose coordinates $(t_1,\dots,t_r)$ on $V_r$,  the additive group of dimensi\'{o}n $r$ is ${\mathbb A}^r=Spec\,k[t_1,\dots,t_r]$.

\begin{subsection}{Tangent space}

 If $G$ is a group with identity element $e$, we denote by $T_eG$ its tangent space at $e$ (Lie algebra). We will use that the elements of $T_eG$ can be calculated with dual numbers: $Der_k({\mathcal O}_G,k)\simeq $

 $\{ f:{{\mathcal O}}_G\to k[\epsilon]=k[x]/x^2 , \text{ $k$-algebra morphism  and }f^{-1}(\epsilon)=\frak m_e \}$.

One can see this in Chap II \S 4 of [DG].

Thus it is proven that  \  $T_eAut_kX\simeq Der_k({\mathcal O}_X,{\mathcal O}_X)$ [MO]. Similarly one can prove that for a $N$-graded algebra $A$,  $T_eAut_{k-gr}A\simeq Der_{N-gra}(A)$ (Definition 3.5).

\begin{prop} Let $H_1$ and $H_2$ be  connected subgroups of an affine group $G$. If \  $T_eH_1\subseteq T_eH_2$ \  inside $T_eG$, then $H_1\subseteq H_2$.\end{prop}
\begin{pf} One can see this in [DG]. \end{pf}
\end{subsection}

\section{Toric variety and divisors}

In this article the field $k$ will be algebraically closed and of characteristic 0 and $(X,{\mathcal O}_X)$ a complete toric variety over $k$ with torus $T_M$. $M$ will be a free ${\mathbb Z}$-module of rank $n$ and $ M^*=Hom_{{\mathbb Z}}(M,{\mathbb Z})$ its dual module.

 If $S\subset M$, we denote \ $k[S]=\{\underset{\alpha\in S}{\sum} \lambda_\alpha x^\alpha \  / \  \lambda_\alpha\in k\text{ and almost all zero }\}$.

  We denote $\Delta$ to the  fan of $X$ whose elements are polyhedral cones of \  $M^*_{\mathbb Q}=M^*\otimes_{\mathbb Z}\mathbb Q$.  This is \  $X=\underset{\sigma\in \Delta}{\bigcup} U_\sigma$ \  where \   $S_\sigma=\sigma^\vee\cap M=$

  $=\{\alpha\in M \, / \, v(\alpha)\geq 0 \text{ for all } v\in\sigma\}$ \  and \  $U_\sigma=Spec\,k[S_\sigma]$. $U_M=Spec\,k[M]$ is an open set contained in each open set $U_\sigma$.

For each ray $\sigma_1$ of the fan there is  only one epimorphism \  $v_i:M\to {\mathbb Z}$ \  generator of the ray; that is \  $\sigma_1=<v_i>_{\mathbb Q^+}$. We will denote the ray $\sigma_1$ by $v_i$ instead of $<v_i>_{\mathbb Q^+}$.  The  rays of $\Delta$ are in correspondence with the  irreducible hypersurfaces of $X$ invariant  by the group $T_M$.

Henceforth \  $\Delta_1=\{v_1,\dots,v_r\}$ \   will be the set of rays of $\Delta$ \  and  \  $H_1=H_{v_1},\dots,H_r=H_{v_r}$ \ the corresponding $T_M$- invariant irreducible hypersurfaces of $X$.
 We have that \ $U_M=X \, \backslash \, (H_{1}\cup\dots\cup H_r)$.

These hypersurfaces $H_i$ define discrete valuations of $k(X)$ (fraction field  of $X$) as follows:
 If \  $f=\underset{\alpha\in M}{\sum} \lambda_\alpha x^\alpha\in k[M]$, $v_i(f)= \underset{\lambda_\alpha\neq 0}{\text{ min }}\{v_i(\alpha) \}$  (*).

As\  $(X,{\mathcal O}_X)$ \  is a complete variety,  $k=\Gamma(X,{\mathcal O}_X)=\underset{\sigma\in \Delta}{\bigcap} k[S_\sigma]=k[\underset{v\in \Delta_1}{\bigcap} S_v]$. That means that $\underset{v_i\in \Delta_1}{\bigcap} S_{v_i}=0$ and so \  $<v_1,\dots,v_r>_{\mathbb Q^+}=M_{\mathbb Q}^*$.

We will use often that if $v_i(\alpha)\geq 0$ for all $v_i\in \Delta_1$, then $\alpha=0$.
\begin{prop} For each ray $v_i\in \Delta_1$, it holds

 $$\Gamma(X,{\mathcal O}_X(H_{i}))=k[P_{v_i}]=\{\underset{\alpha\in P_{v_i}}{\sum} \lambda_\alpha x^\alpha \  / \  \lambda_\alpha\in k\}$$

where \  $P_{v_i}=\{\alpha\in M \ / \ {v_i}(\alpha)=-1 \text{ and } v_j(\alpha)\geq 0 \text{ for all   } v_j\in \Delta_1 \text{ distinct from } {v_i}\}$
\end{prop}

\begin{pf} $ \Gamma(X,{\mathcal O}_X(H_{v_i}))\subset \Gamma(U_M,{\mathcal O}_X(H_{v_i}))=k[M]$. $f\in \Gamma(X,{\mathcal O}_X(H_{v_i}))$ \  if and only if \  $f\in k[M]$ \  and \  $v_i(f)\geq -1$  and $v_j(f)\geq 0 $ for all $v_j\in \Delta_1\text{ different from } v_i$. By \  (*) \   $f\in P_{v_i}$.
\end{pf}

 In [F] is proved that the  Weil divisor group of $X$ modulo  linear equivalence, $N=A(X) $, is generated by the classes of the hypersurfaces \  $[H_{1}],[H_{2}],\dots,[H_{r}]$. Indeed the following sequence is exact:
$$0\to M\overset{v}{\to} {{\mathbb Z}}^r\overset{\pi}{\to} N\to 0$$

where \  $v(\alpha)=(v_1(\alpha),\dots,v_r(\alpha))$ \  and \ $\pi(n_1,\dots,n_r)=[\overset{r}{\underset{i=1}{\sum}} n_iH_{i}]$.

 We denote $Pic(Z)$ the group of invertible sheaf of $Z$ modulo isomorphism. If $Z$ is nonsingular, then $Pic(Z)$ is $A(Z)$.

We will call  $U$ the open set of nonsingular points  $X$. As $U$ contains every point of   codimension 1 of $X$, one has that \  $Pic(U)=A(U)=A(X)=N$ and for each divisor $D$ of $X$, $\Gamma(X,{\mathcal O}_X(D))=\Gamma(U,{\mathcal O}_U(D))$.

From now on, we denote \  $Z\times_{Spec \, k}  Z'=Z\times Z'$ \  for each $k$-scheme $Z$ and $Z'$.

 \begin{prop} Let  $Z$ be  a smooth variety of fraction field  $k(Z)$ and

 $\pi:U\times Z\to U$ ,
 \  $p:U\times Z\to Z$ \ the canonical projections.  The morphism   $Pic(U)\times Pic(Z)\to Pic(U\times Z)$ \  given by:
 $(L_1,L_2)\to \pi^*L_1\otimes_{{\mathcal O}_{U\times Z}} p^*L_2$ \  is an isomorphism.
 \end{prop}

 \begin{pf} We denote  $\bar{H}_i=\pi^{-1}H_i$ ,
 $\pi^{-1}(U_M)=\bar{U}_M$ the open set of  $X\times Z$ \  whose complement set is \  $\bar{H}_1\cup\dots\cup \bar{H}_r$.

Let us see what discrete valuations of  $k(U\times Z)$ \  center  at  codimension 1:

 Let $H$ be an integral hypersurface   of $U\times Z$ \  defining the discrete valuation \  ${\mathcal O}_{v_H}$.

   Computing  dimension, either  $p(H)=Z$   or \  $p(H)=H'$ \   is a hypersurface of  $Z$. In the latter case  \ $H=p^{-1}(H')$ \  and \   ${\mathcal O}_{U\times Z}(H)=p^*{\mathcal O}_{ Z}(H')$.

 Therefore either \  $v_H$ \ centers  on  \  $U\times  Spec\,k(Z)$ \   or \  ${v_H}=v_{p^{-1}H'}$ \  for some hypersurface $H'\subset Z$.

This morphism is surjective: Let  $H\neq p^{-1}H'$ be an integral hypersurface of $U\times  Z$. On the toric variety over the field $k(Z)$, $U\times  Spec\,k(Z)$, \  is   $H\equiv \sum n_i\bar{H}_i$. Therefore on  $U\times  Spec\,k(Z)$,  $H- \sum n_i\bar{H}_i=$ div\,$f$ for some $f\in k(X\times  Z)$  \  and so \   div\,$f=H- \sum n_i\bar{H}_i+\sum m_jp^{-1}H'_j$  \  on $X\times  Z$. We can conclude that  $[H]=\sum n_i[\bar{H}_i]+\sum m_j[p^{-1}H'_j]$.

 The morphism is injective: If \  $\pi^*L_1\otimes p^*L_2\simeq {\mathcal O}_{U\times  Z}$,  taking a rational point $z$ of $Z$ and restricting to $p^{-1}(z)=U\times  \{z\}$, one has that $L_1\simeq {\mathcal O}_U$. Analogously  taking a rational point of $U$,  $L_2\simeq {\mathcal O}_{Z}$.

\end{pf}
\begin{cor} If \ $\tau:X\to X$ \  is an automorphism in \  $Aut^0_kX$, then for each $i$ \  $\tau(H_i)$ is linearly equivalent to $H_i$.
\end{cor}
\begin{pf}Let \  $Z=Aut_k^0X$ \  and \  $\phi: X\times Z\to X\times Z$ the universal automorphism; this is \  $\phi(x,\tau)=(\tau(x),\tau)$. The  non singular locus of \ $X\times Z$ \  is \  $U\times Z$ \  and so $\phi$ is an automorphism of \  $U\times Z$. If \  $L_i={\mathcal O}_U(H_i)$, then \  $\tau^*L_i={\mathcal O}_U(\tau^{-1}H_i)$. We only have to prove that \  $\tau^*L_i\simeq L_i$. By the above proposition \  $\phi^*\pi^*L_i\simeq p^*L\otimes\pi^*L'$. If we restrict this isomorphism to $U\times\{Id\}$, we have that \  $L_i\simeq {\mathcal O}_U\otimes L'=L'$ \  and if we restrict  to $U\times\{\tau\}$, we have that \  $\tau^*L_i\simeq {\mathcal O}_U\otimes L'=L'$.
\end{pf}
\begin{prop} With the  hypothesis of proposition 2.2, $\Gamma(U\times Z,{\mathcal O}_{U\times Z})^\times=\Gamma( Z,{\mathcal O}_{ Z})^\times$.\end{prop}
\begin{pf} We have a birational morphism $U_M\times Spec\,k(Z)=Spec\,k(Z)[M] \to U\times Z$. Therefore \  $f\in k(Z)[M]^\times$. Hence $f=\lambda_\alpha x^\alpha $ where $\lambda_\alpha\in k(Z)$ and $\alpha\in M$. Besides $v_H(f)=0$ for each valuation $v_H$ defined by an irreducible  hypersurface $H\subset U\times Z$. If \  $H=H_i\times Z$, $v_H(f)=v_i(\alpha)=0$. Therefore $\alpha=0$. If $H=U\times \bar{H}$, $v_H(f)=v_{\bar{H}}(\lambda_\alpha)=0$ for every irreducible hypersurface $\bar{H}\subset Z$. We conclude that \  $f=\lambda_\alpha\in \Gamma(Z,{\mathcal O}_{Z})^\times$.\end{pf}

\section{Derivations}

If $X$ is a $k$-scheme, we denote \   $Der_k({\mathcal O}_X)=Der_k({\mathcal O}_X,{\mathcal O}_X)$ \  the $k$-vectorial space of the $k$-derivations from ${\mathcal O}_X$ to ${\mathcal O}_X$.

If \  $w\in M^*\otimes_{{\mathbb Z}} k=Hom_{{\mathbb Z}}(M,k)$, this defines a  derivation of $k[M]$, $D_w$ , given by \  $D_wx^\alpha= w(\alpha)x^\alpha$ \  for each $\alpha\in M$. This derivation is  a derivation of $k[S]$ for each semigroup $S\subset M$ \  and so  this is a global derivation of $X$ for any toric variety $X$ of group $T_M$.

\begin{prop} a) The morphism \  $M^*\otimes_{{\mathbb Z}}k[M] \to Der_k(k[M],k[M])$ \  given by \  $v\otimes f\to fD_v$ \  is  an isomorphism.

b) Let \  $v_i\in M^*$ ,   $S_{v_i}=\{\alpha\in M \ / \ v_i(\alpha)\geq 0\}$ \ and \ $P_{v_i}=\{\alpha\in M \ / v_i(\alpha)= -1\}$.

 $$Der_k(k[S_{v_i}],k[S_{v_i}])=M^*\otimes_{{\mathbb Z}}k[S_{v_i}] \oplus k[P_{v_i}]\cdot D_{v_i}$$

c) Let $\sigma$ be  a polyhedral cone generated by  $v_1,\dots ,v_s\in M^*$, $S_\sigma=M\cap \sigma^\vee$   and \

 $ P_{v_i}=\{\alpha\in M \, / \ v_j(\alpha)\geq 0 \text{ for all } j\neq i \text{ and } v_i(\alpha)= -1\}$. We have that

  $$Der_k(k[S_\sigma],k[S_\sigma])=M^*\otimes_{{\mathbb Z}}k[S_\sigma]\oplus k[P_{v_1}]\cdot D_{v_1}\oplus\dots\oplus k[P_{v_s}]\cdot D_{v_s}$$
\end{prop}
\begin{pf}a)  If $e_1,\dots,e_n$ is a basis of $M$, $w_1,\dots,w_n$ its dual basis and $x_i=x^{e_i}$, then $k[M]=k[x_1,\dots,x_n,\frac{1}{x_1\cdot{ \dots }\cdot x_n}]$ \  and  \  $\{D_{w_i}= x_i\frac{\partial}{ \partial x_i}\}_{1\leq i\leq n}$ \ is a basis of derivations of $k[M]$.

b) As $k[M]$ is a localization of $k[S_{v_i}]$, we have that  \  $Der_k(k[S_{v_i}],k[S_{v_i}])\subset$  $Der_k(k[M])$.
 A derivation $D$ of $k[M]$ is a derivation of $k[S_{v_i}]$ if \ $Dx^\alpha \in k[S_{v_i}]$ \  for all $\alpha\in S_{v_i}$.

 Let us  see when \  $D=x^\beta D_w$ \  is a derivation of $k[S_{v_i}]$:

 Let $\alpha\in S_{v_i}$; that is $v_i(\alpha)\geq 0$. $Dx^\alpha=x^\beta D_wx^\alpha=x^{\beta+\alpha} w(\alpha)\in k[S_{v_i}]$ \  if and only if \  $w(\alpha)=0$ \  or \  $v_i(\beta+\alpha)=v_i(\beta)+v_i(\alpha)\geq 0$.

  If $w(\alpha)=0$ for all $\alpha\in M$ such that $v_i(\alpha)=0$,  then $w=\lambda v_i$ for some $\lambda \in k$ and \  $D=\lambda x^\beta D_{v_i}$. Besides, if  $v_i(\alpha)= 1$, $Dx^\alpha\in k[S_{v_i}]$ if and only if $v_i(\beta)\geq -v_i(\alpha)= -1$ and so $\beta \in P_{v_i}$ or $\beta\in S_{v_i}$.

 If there exists $\alpha\in M$ such that $v_i(\alpha)=0$ and $w(\alpha)\neq 0$, $v_i(\beta)+v_i(\alpha)=v_i(\beta)\geq 0$ and so \   $D\in M^*\otimes_{{\mathbb Z}} k[S_{v_i}]$.

 c)  We consider $Der_k(k[S_\sigma])$ and $Der_k(k[S_{v_i}])$ as subsets of $Der_k(k[M])$.
 A derivation $D$ of $k[M]$ is a derivation of $k[S_\sigma]$ if \ $Dx^\alpha \in k[S_\sigma]=\overset{s}{\underset{i=1}{\bigcap}} k[S_{v_i}]$ \  for all $\alpha\in S_\sigma$.
 That is,  $Der_k(k[S_\sigma],k[S_\sigma])=\overset{r}{\underset{i=1}{\bigcap}} Der_k(k[S_{v_i}],k[S_{v_i}])$ and we can conclude by b).

\end{pf}

\begin{thm} Let $(X,{\mathcal O}_X)$ be a  toric variety  with rays \   $v_1,\dots,v_r$  \  and

$ P_{v_i}=\{\alpha\in M \, / \ v_j(\alpha)\geq 0 \text{ for all } j\neq i \text{ and } v_i(\alpha)= -1\}$.
 $$Der_k({\mathcal O}_X,{\mathcal O}_X)=M^*\otimes_{{\mathbb Z}} k\oplus k[P_{v_1}]\cdot D_{v_{1}}\oplus\dots\oplus  k[P_{v_r}]\cdot D_{v_{r}}$$
\end{thm}
\begin{pf} If \ $D:{\mathcal O}_X\to{\mathcal O}_X$ \ is a derivation, taking sections over $U_M$, we have that $D$ is a derivation of $k[M]$. A derivation $D$ of $k[M]$ is a global derivation of $X$ if and only if $D(k[S_\sigma])\subset k[S_\sigma]$ \  for all \  $\sigma\in \Delta$. As each $k[S_{v_i}]$ is a localization of $k[S_\sigma]$ \ and \ $S_\sigma=\underset{v_i\in\sigma}\cap S_{v_i}$, it is enough to verify that \ $D\in Der_k(k[S_{v_i}])$ . If \  $D=x^\beta\cdot D_w$ \  and \  $x^\beta\notin k[S_{v_i}]$, then  by the above proposition, we conclude  that \  $w=v_i$ for some $i$ \  and \  $\beta\in P_{v_i}$.
\end{pf}
\begin{cor} $T_M$ is a maximal torus of $Aut_k^0X$.
\end{cor}
\begin{pf}If $T_M\subset T'$ and $T'$ is a torus, then $T_M$ acts trivially on $T_eT'$ and so $T_eT'=M^*\otimes_{{\mathbb Z}}k=T_eT_M$.
\end{pf}

\begin{dfn} We call \underline{root system} of \  $Aut^0_kX$ \  to the set \  $\overset{r}{\underset{i=1}{\bigcup}}  P_{v_i}$.
\end{dfn}
\begin{subsection}{Graded derivations}

\begin{dfn}Let $N$ be a ${\mathbb Z}$-module and $A=\underset{n\in N}{\oplus} A_n$  a $N$-graded $k$-algebra. We will say that a $k$-derivation  \   $D:A\to A$ is graded if \ $D(A_n)\subseteq A_n$ for all $n\in N$.

We will denote it by \  $Der_{N-gra}(A)$.
\end{dfn}
\begin{thm} Let $\sigma$ be  a polyhedral cone generated by  $v_1,\dots ,v_s\in L^*$, $S_\sigma=L\cap \sigma^\vee$ \  and \  $P_{v_i}$ \  as proposition 3.1 c). Let \   $deg:L\to N$ \   be a morphism  and we consider $k[S_\sigma]$ as a $N$-graded algebra via $deg$; that is \  $deg(x^s)=deg(s)$ for all $s\in S_\sigma$. If \  $Ker\ deg=N'$, it holds that

  $ Der_{N-gra}(k[S_\sigma])=L^*\otimes_{{\mathbb Z}} k[S_\sigma\cap N']\oplus k[P_{v_1}\cap N']\cdot D_{v_1}\oplus \dots\oplus k[P_{v_s}\cap N']\cdot D_{v_s}$
\end{thm}
\begin{pf} By the proposition 3.1 c)

$ Der_{k}(k[S_\sigma])=L^*\otimes_{{\mathbb Z}}k[S_\sigma]\oplus  k[P_{v_1}]\cdot D_{v_1}\oplus \dots\oplus k[P_{v_s}]\cdot D_{v_s}$.

We can conclude taking into account that for each $w\in L^*$ and  $ \beta\in L$, $x^\beta D_w$ \ is a $N$-graded derivation if and only if \  $\beta\in Ker\, deg$.
\end{pf}

\begin{cor} Let \  $A=k[z_1,\dots,z_l]$ \  a polynomial ring.  We consider $A$ as a ${\mathbb Z}$-graded algebra where the degree of each variable $z_i$ is $1$. $Der_{N-gr}A$ is the $k$-vectorial space generated (inside $Der_kA$) by the derivations $\{z_i\frac{\partial}{\partial z_j}\}_{1\leq i,j\leq l}$.
\end{cor}
\begin{pf}  By above theorem we have that \  $Der_{{\mathbb Z}-gra}A=({\mathbb Z}^l)^*\otimes_{{\mathbb Z}}k \, \oplus <z^{r_{ij}}D_{v_j}>$  \  where \  $v_j:{\mathbb Z}^l\to{\mathbb Z}$ \  is the $j$-th canonical projection and \  $r_{ij}\in {\mathbb Z}^l$ \  is an element of degree $0$ such that \  $v_j(r_{ij})=-1$ and $v_k(r_{ij})\geq 0$ for all $k\neq j$. We only  have to check that for $i\neq j$, $z_i\frac{\partial}{\partial z_j}=z^{r_{ij}}D_{v_j}$ \  where \ $z^{r_{ij}}=\frac{z_i}{z_j}$ \  and for $i=j$, $z_j\frac{\partial}{\partial z_j}=D_{v_j}$  generate the derivations  corresponding to ${{\mathbb Z}^l}^*$.

\end{pf}

\end{subsection}

\section{Cox ring}{}
We carry on with notations of above sections. Let $n\in {{\mathbb Z}}^r$. We say that $n\geq 0$ if $n\in{\mathbb N}^r$. We recall that we have the exact sequence:
   \  \   $0\to M\overset{v}{\to} {{\mathbb Z}}^r\overset{\pi}{\to} N\to 0$

and  for each $\alpha\in M$,   $v(\alpha)\geq 0$ if and only if $\alpha=0$.

\begin{dfn} We call Cox sheaf ring to ${{\mathcal A}}=\underset{(n_1,\dots,n_r)\in{{\mathbb Z}}^r}{\bigoplus} {\mathcal O}_X(n_1H_1+\dots+n_rH_r)$

and Cox general ring of $X$ to $A=\Gamma(X,{\mathcal A})$.

 \end{dfn}

 \begin{prop} $A\simeq k[S]$   \  where \   $S=\{(\alpha,n)\in M\oplus {\mathbb Z}^r  \  / \  v(\alpha)+n\geq 0\}$.
 \end{prop}

 \begin{pf} $\Gamma(X,{\mathcal A})\subset \Gamma(U_M,{\mathcal A})=\underset{m\in{\mathbb Z}^r}{\bigoplus} k[M]_m$ \ where $k[M]_m=k[M]$.

 Therefore \  $A\subset k[M][{\mathbb Z}^r]=k[M\oplus {\mathbb Z}^r]$. We just have to compute the elements of $A$ of degree $m\in{\mathbb Z}^r$:
  For each $n=(n_1,\dots,n_r)\in{\mathbb Z}^r$,
  $f\in \Gamma(X,{\mathcal O}_X(n_1H_{1}+\dots +n_rH_{r}))$  \  if and only if \  $f\in k[M]$ \  and \  $v_i(f)\geq -n_i$ \   for all $1\leq i\leq r$.

 By \  (*) \   $f=\underset{\alpha,v(\alpha)+n\geq 0}{\sum} \lambda_\alpha x^\alpha$.
 \end{pf}

  If we consider ${\mathbb N}^r$ as   $0\oplus {\mathbb N}^r\subset S$, we have that  $A_C=k[{\mathbb N}^r]\subset k[S]$. This subring is the ring constructed by Cox in [C].

  We will consider $A=k[S]$ as ${\mathbb Z}^r$-graded algebra and $A_C$ as $N$-graded algebra via the morphism $\pi:{\mathbb Z}^r\to N$. Its groups of graded automorphisms will be denote  respectively by $Aut_gA$ and $Aut_gA_C$.

   We denote the elements of $k[M]$ by $\underset{\alpha\in M}{\sum} \lambda_\alpha x^\alpha$ and the elements of $k[{\mathbb Z}^r]$ by $\underset{m\in {\mathbb Z}^r}{\sum} \mu_m y^m$. Thus \ $A_C=k[y_1,\dots,y_r]$ and the elements of $k[S]$ are denoted by \  $\underset{(\alpha,m)\in S}{\sum} \lambda_{\alpha,m}x^\alpha y^m$.
 We can recover the fraction field of $X$, $k(X)$, with $A$ and $A_C$:

 \begin{prop} If $B$ is a graded algebra, we denote

 $k(B)_0=\{\dfrac{p}{q} \  / \ p, q \text{ are homogeneous elements of $B$ with the same degree} \}$.

 a) $k(A)_0=k(X)$.

b)   $k(X)\simeq k(A_C)_0$ \  identifying \  $x^\alpha$ \  with  \  $y^{v(\alpha)}$ \ for each $\alpha\in M$.
\end{prop}
\begin{pf}

 a) If \  $p,q$ are homogeneous elements of $A$ of degree \   $n\in{\mathbb Z}^r$, then

$\dfrac{p}{q}= \dfrac{\underset{\alpha\in M}{\sum} \lambda_{\alpha} x^\alpha y^n}{\underset{\alpha\in M}{\sum} \mu_{\alpha} x^{\alpha}y^n}=\dfrac{\underset{\alpha\in M}{\sum}\lambda_{\alpha} x^{\alpha}}{\underset{\alpha\in M}{\sum} \mu_{\alpha} x^{\alpha}}\in k(X)$.

Conversely, we only have to prove that \  $M\subset S-S=\{s_1-s_2 \   / \  s_1,s_2\in S\}$. $0\oplus {\mathbb N}^r\subset S$ and so $0\oplus {\mathbb Z}^r\subset S-S$. If $\alpha\in M$, $\alpha=(\alpha,-v(\alpha))+(0,v(\alpha))\in S-S$.

b) $v(\alpha)=n_1-n_2$ \  where $n_1,n_2\in {\mathbb N}^r$. So $y^{v(\alpha)}=\dfrac{y^{n_1}}{y^{n_2}}\in k(A_C)_0$. If \  $\dfrac{p}{q}$ \   is a fraction such that \ $p$ and $q$ \ are homogeneous elements of $A_C$ with the same degree $\pi(n)\in N$, then \ $\dfrac{p}{q}=\dfrac{\underset{\alpha\in M}{\sum} \lambda_{\alpha} y^{v(\alpha)+n}}{\underset{\alpha\in M}{\sum} \mu_{\alpha} y^{v(\alpha)+n}}=\dfrac{\underset{\alpha\in M}{\sum}\lambda_{\alpha} y^{v(\alpha)}}{\underset{\alpha\in M}{\sum} \mu_{\alpha} y^{v(\alpha)}}$ \  \  is identified with  \
 \  $\dfrac{\underset{\alpha\in M}{\sum}\lambda_{\alpha} x^{\alpha}}{\underset{\alpha\in M}{\sum} \mu_{\alpha} x^{\alpha}}\in k(X)$.

\end{pf}

  \begin{lem} Let $K$ be a $k$-algebra.

  a) $K[S]$ \  is generated by the elements of the form  \  $x^\alpha y^{-v(\alpha)}\cdot y^{m}$ \  where \  $\alpha\in M$ \ and \ $m\in {\mathbb N}^r$.

  b) If \  $\tau:K[S]\to K[S]$ \  is   a ${\mathbb Z}^r$-graded automorphism, then
  for each $\alpha\in M$, $\tau(x^\alpha y^{-v(\alpha)})=\lambda_\alpha \, x^\alpha y^{-v(\alpha)}$ \  for some $\lambda_\alpha\in K$.

   c) For each \  $m\in{\mathbb N}^r$, we denote \  $Y_m=\{\beta\in M \   / \ v(\beta)+m\geq 0\}$.

   If \  $\tau:K[{\mathbb N}^r]\to K[{\mathbb N}^r]$ \  is   a $N$-graded automorphism, then \ $\tau(y^m)=\underset{\beta\in Y_m}{\sum} \lambda_{\beta} y^{m+v(\beta)}$ \  for each \  $m\in{\mathbb N}^r$.

  \end{lem}
  \begin{pf} a) $K[S]$ \  is  generated by the elements of the form \  $x^\alpha\cdot y^n$ where $\alpha\in M$, $n\in {\mathbb Z}^r$ and $v(\alpha)+n\geq 0$.
 $x^\alpha\cdot y^n=x^\alpha y^{-v(\alpha)}\cdot y^{n+v(\alpha)}=x^\alpha y^{-v(\alpha)}\cdot y^{m}$ \  where \  $m\in {\mathbb N}^r$.

  b) $\tau(x^\alpha y^{-v(\alpha)})$ \  is an element of degree \  $-v(\alpha)$.

  $x^\beta\cdot y^n$ has degree $-v(\alpha)$ if $n=-v(\alpha)$. But \  $0\leq v(\beta)+n=v(\beta)-v(\alpha)=v(\beta-\alpha)$. Therefore \  $\beta-\alpha=0$ and $x^\beta\cdot y^n=x^\alpha y^{-v(\alpha)}$.

  c) $\tau(y^m)=\sum \lambda_n y^n$ \  where $y^n$ and $y^m$ have the same $N$-degree. Therefore \  $n=m+v(\beta)$ \ for some $\beta\in M$. As $n\geq 0$, one concludes that $\beta\in Y_m$.

  \end{pf}

  Let  \  $\tau\in Aut_{N-gr}K[{\mathbb N}^r]$ \  such that \  $\tau(y^m)=\underset{\beta\in Y_m}{\sum} \lambda_{m,\beta} y^{m+v(\beta)}$.

  We denote \  $i(\tau)=\bar{\tau}:K[S]\to K[S] $  \  the ${\mathbb Z}^r$-graded morphism given by \
   $\bar{\tau}(x^\alpha y^{-v(\alpha)})=x^\alpha y^{-v(\alpha)}$ for each $\alpha\in M$ \  and  \  $\bar{\tau}(y^m)=\underset{\beta\in Y_m}{\sum} \lambda_{m,\beta}x^{\beta} y^{m}$ \ for each \  $m\in {\mathbb N}^r$.

   This is straightforward to prove that $\bar{\tau}$ is a morphism; that is \  $\bar{\tau}(x^\alpha\cdot y^n \cdot x^\beta y^m)=\bar{\tau}(x^\alpha\cdot y^n)\cdot \bar{\tau}(x^\beta\cdot y^n)$ \  for all \  $(\alpha,m),(\beta,m)\in S$.

   \begin{lem} a) $\overline{\tau_1\circ\tau_2}=\bar{\tau}_1\circ\bar{\tau}_2$ \ for each $\tau_1,\tau_2\in Aut_gA_C $.

   b) $\bar{\tau}$ is an automorphism.
   \end{lem}

  \begin{pf} a) If \  $\tau_1(y^n)=\underset{\alpha\in Y_n}{\sum} \lambda_{n,\alpha} y^{n+v(\alpha)}$ \  and \  $\tau_2(y^n)=\underset{\alpha\in Y_n}{\sum} \mu_{n,\alpha} y^{n+v(\alpha)}$, then

\medskip

$\tau_1\circ \tau_2(y^n)=\tau_1(\underset{\alpha\in Y_n}{\sum} \mu_{n,\alpha} y^{n+v(\alpha)})=\underset{\alpha\in Y_n}{\sum} \mu_{n,\alpha} \underset{\beta\in Y_{n+v(\alpha)}}{\sum} \lambda_{n+v(\alpha),\beta}  y^{n+v(\alpha)+v(\beta)}$

\medskip

$=\underset{\alpha\in Y_n,\beta\in Y_{n+v(\alpha)}}{\sum} \mu_{n,\alpha} \lambda_{n+v(\alpha),\beta}  y^{n+v(\alpha+\beta)}$.

  \medskip

  Therefore \  $\overline{\tau_1\circ\tau_2}(y^n)=\underset{\alpha\in Y_n,\beta\in Y_{n+v(\alpha)}}{\sum} \mu_{n,\alpha}  \lambda_{n+v(\alpha),\beta} x^{\alpha+\beta} y^{n}$.

  \medskip

  $\bar{\tau}_1\circ \bar{\tau}_2(y^n)=\bar{\tau}_1(\underset{\alpha\in Y_n}{\sum} \mu_{n,\alpha} x^\alpha y^{n})=\underset{\alpha\in Y_n}{\sum} \mu_{n,\alpha} x^\alpha y^{-v(\alpha)} \bar{\tau}_1(y^{n+v(\alpha)})=$

  \medskip

  $\underset{\alpha\in Y_n}{\sum} \mu_{n,\alpha} x^\alpha y^{-v(\alpha)}\underset{\beta\in Y_{n+v(\alpha)}}{\sum}  \lambda_{n+v(\alpha),\beta} x^{\beta}  y^{n+v(\alpha)}=\underset{\alpha\in Y_n,\beta\in Y_{n+v(\alpha)}}{\sum} \mu_{n,\alpha}  \lambda_{n+v(\alpha),\beta} x^{\alpha+\beta}y^{n}$.

  b) $\overline{Id}=Id$ \  and so \  $Id=\overline{\tau\circ\tau^{-1}}=\bar{\tau}\circ\overline{\tau^{-1}}$.

  \end{pf}

  \begin{thm} For each  \   $\tau\in Aut_gA$, we denote by \   $p(\tau)$ the element of $ T_M$ \   such that \   $p(\tau)(\alpha)=\lambda_\alpha$  \    if \    $\tau(x^\alpha y^{-v(\alpha)})=\lambda_\alpha x^\alpha y^{-v(\alpha)}$.

  $0\to Aut_gA_C \overset{i}{\to}Aut_gA\overset{p}{\to}T_M \to 0$ \  is a split exact sequence.
 \end{thm}
  \begin{pf} Let \  $\bar{\tau}\in Ker\, p$. For each $n\in{\mathbb N}^r$, $\bar{\tau}(y^n)$ has degree $n$.
 Therefore \

$\bar{\tau}(y^n)=\underset{\alpha\in Y_n}{\sum} \lambda_{n,\alpha} x^\alpha y^{n}$ \  and   $\tau\in Aut_gA_C $ \   defined by  \  $\tau(y^{n})=\underset{\alpha\in Y_n}{\sum} \lambda_{n,\alpha} y^{n+v(\alpha)}$ \ is such   that \   $i(\tau)=\bar{\tau}$.

The section of $p$ is the action of $T_M$ on $k[S]\subset k[M\oplus {\mathbb Z}^r]$ as $M$-graded algebra.
\end{pf}
\begin{cor} $Aut_gA_C  \, / \, T_{N}\simeq Aut_gA/T_{{\mathbb Z}^r}$ \  and \
 $Aut^0_gA_C  \, / \, T_{N}\simeq Aut^0_gA/T_{{\mathbb Z}^r}$

\end{cor}
\begin{pf}  We have the exact sequences

 $$\begin{matrix}0\to&T_{N}&\to &T_{{\mathbb Z}^r}&\to &T_M\to 0                    \\
 &\bigcap & & \bigcap&& \hskip -1cm ||\\
 0\to &Aut_gA_C &\overset{i}{\to}&Aut_gA&\overset{p}{\to}&T_M \to 0\end{matrix}$$

and taking quotient groups  we conclude.

 Because of being $T_{{\mathbb Z}^r}$ connected, this is included in $Aut^0_gA\simeq Aut_g^0A_C\rtimes T_M$ and  we have the similar exact sequence as above.
\end{pf}

\begin{subsection}{Tangent space of \  $Aut_gA_C \, / \, T_N$}
With the previous notations, $A_C=k[{\mathbb N}^r]=k[y_1,\dots,y_r]$ ; $v_i:{\mathbb Z}^r\to {\mathbb Z}$ is the canonical projections; that is $v_i(n_1,\dots,n_r)=n_i$. If $\sigma$ is polyhedral cone generated by $\{v_1,\dots,v_r\}$, we have that ${\mathbb N}^r={\mathbb Z}^r\cap\sigma^\vee$ and $Ker\,\pi\cap {\mathbb N}^r=0$.

By the theorem 3.6 the $N$-graded derivations of  $A_C$ are:
$$Der_{g}(A_C)=({\mathbb Z}^r)^*\otimes_{{\mathbb Z}}k\ \oplus \  <y^{r_{ij}}D_{v_i}>$$

 where \  $r_{ij}\in {\mathbb Z}^r$ \  hold that $\pi(r_{ij})=0$ , $v_i(r_{ij})=-1$ \ and \  $v_j(r_{ij})\geq 0$ for all $j\neq i$. A maximal torus of $Aut_g^0A_C$ is $T_{{\mathbb Z}^r}$ and these $r_{ij}$ are \underline{the root system} of $Aut^0_gA_C$.

\begin{rem} $Ker\,\pi=v(M)$ and $v$ gives an one to one correspondence  between the root system  of $Aut^0_kX$ and the root system of $Aut^0_gA_C$. This is: $\alpha\in M$ is a root of $Aut_k^0X$ if and only if $v(\alpha)$ is a root of $Aut_gA_C$.
\end{rem}

We have that \ $T_{{\mathbb Z}^r}/T_{N}\simeq T_M$ \ and so \  $T_e(Aut_gA_C/T_N)=M^*\otimes_{{\mathbb Z}}k \ \oplus \  <y^{r_{ij}}D_{v_i}>$
\end{subsection}

\section{Relation between the groups.}

In this section  we use the notations of section 2.  $U$ is the open set of nonsingular points of $X$ \ and for each $k$-scheme $Y$, $\pi:U\times Y\to U$ \ and \  $p:U\times Y\to Y$ \ are the canonical projections.

Let $L$ be an invertible sheaf of $U$. We denote \  $Aut_LU$ \  to the following functor:

 For each $k$-scheme $Y$, $Aut_LU\,(Y)$ \ is the set of pairs \   $(\tau,\phi)$ \  where

$\tau:U\times Y\to U\times Y$ \  is an automorphism over $Y$ ($p\circ\tau=p$) \  and \  $\phi:\tau^*\pi^*L\simeq \pi^*L$ \  is an isomorphism.

 If \  $f:Y'\to Y$ \  is a  morphism and $\bar{f}=Id\times f$, we have that \   $\pi\circ\bar{f}=\pi$ \  and    the following  diagrams are commutative:

$$  \hskip 2cm \begin{matrix}U\times  Y'&\overset{\bar{f}}{\to}&U\times  Y&& \\
\hskip -0.3cm \bar{\tau} \, \downarrow& &\hskip -0.3cm \tau \, \downarrow& & \\
U\times  Y'&\overset{\bar{f}}{\to}&U\times  Y&\overset{\pi}{\to}& U\\
\hskip -0.3cm p \, \downarrow& &\hskip -0.3cm p \, \downarrow&(1) & \downarrow\\
Y'&\overset{f}{\to}& Y&{\to}& Spec\, k\end{matrix} $$

  $f^*\tau=\bar{\tau}:U\times  Y'\to U\times  Y'$ \ is the automorphism obtained  taking $\times_YY'$ on \  $\tau:U\times  Y\to U\times  Y$ \  and \    $f^*\phi=\bar{\phi}$ \  is the composition of the isomorphisms \  $\bar{\tau}^*\pi^*L=\bar{\tau}^*\bar{f}^*\pi^*L=\bar{f}^*\tau^*\pi^*L\overset{\phi}{\simeq}\bar{f}^*\pi^*L=\pi^*L $.

 This functor is a functor of groups : If \  $(\tau_1,\phi_1),(\tau_2,\phi_2)\in Aut_LU\,(Y)$ , then

  $(\tau_1,\phi_1)\circ (\tau_2,\phi_2)=(\tau_1\circ\tau_2,\phi_1\circ\phi_2)$ \  where \  $\phi_1\circ\phi_2$ \  is the composition of the isomorphisms \  $\tau_2^*\tau_1^*\pi^*L\overset{\phi_1}{\simeq}\tau_2^*\pi^*L\overset{\phi_2}{\simeq} \pi^*L $.

\begin{prop} Let \  $\mathcal N$ \   a  coherent ${\mathcal O}_U$-module.

a) $\Gamma(U,{\mathcal N})\otimes_k{\mathcal O}_Y\simeq p_*\pi^*{\mathcal N}$.

b) The natural morphism \  $f^*p_*\pi^*{\mathcal N} \to p_*\bar{f}^*\pi^*{\mathcal N}$ \  is an isomorphism.

\end{prop}
\begin{pf} a)  Applying  to the  commutative diagram (1) the proposition 9.3 of [H].

b) As \  $\pi\circ\bar{f}=\pi$, one has that  \  $p_*\bar{f}^*\pi^*{\mathcal N}=p_*\pi^*{\mathcal N}\simeq \Gamma(U,{\mathcal N})\otimes_k{\mathcal O}_{Y'}$ \ applying a). By the other hand \  $f^*p_*\pi^*{\mathcal N} \simeq f^*(\Gamma(U,{\mathcal N})\otimes_k{\mathcal O}_Y)=\Gamma(U,{\mathcal N})\otimes_k{\mathcal O}_{Y'}$.

\end{pf}

  Let $V=\Gamma(U,L)$. We want to show that there exists a group morphism between the  functors $Aut_LU$ and $Aut_kV$:

  Given  $(\tau,\phi)\in Aut_LU \,(Y)$ \ we define \ the automorphism \

 $\tau_{\phi}: p_*\pi^*L=V\otimes_k{\mathcal O}_{Y }\to V\otimes_k{\mathcal O}_{Y }=p_*\pi^*L$ \  as composition of the isomorphisms:
 \  $  p_*\pi^*L\overset{\tau}{\simeq}  p_*\tau_*\tau^*\pi^*L=p_*\tau^*\pi^*L\overset{\phi }{\simeq} p_*\pi^*L$.

\begin{thm} a) The correspondence \  $(\tau,\phi)\to \tau_\phi$ \  is functorial. That is, if  $(\tau,\phi)$ is as above and \  $f:Y'\to Y$ is a morphism,  \  $f^*(\tau_\phi)=(f^*\tau)_{f^*\phi}$.

b) The previous correspondence is a group morphism:

   If \  $\tau_1,\tau_2\in Aut_Y(U\times  Y)$ \  and \  $\phi_1:\tau_1^{*}\pi^* L\to \pi^* L \ , \ \phi_2:\tau_2^{*}\pi^* L\to \pi^* L $ \  are isomorphisms , then
 \  $(\tau_2\circ\tau_1)_{\phi_2\circ \phi_1}=\tau_{2\phi_2}\circ \tau_{1\phi_1}$.

\end{thm}

\begin{pf} a) With the previous notations,
 the  morphism  $f^*\tau_\phi$ is obtained taking $f^*$ on
 \  $  p_*\pi^*L\overset{\tau}{\simeq}  p_*\tau_*\tau^*\pi^*L=p_*\tau^*\pi^*L\overset{\phi }{\simeq} p_*\pi^*L$.

   One has that  $ p\circ \bar{f}=f\circ p$ , $ \pi\circ \bar{f}=\pi$, $\bar{f}\circ\bar{\tau}=\tau\circ\bar{f}$ and one concludes from the following commutative diagrams:

 $$ \begin{matrix}
 V\otimes_k{\mathcal O}_{Y'}=&f^* p_*\pi^*L&\overset{\tau}{\simeq}&f^*p_*\tau^*\pi^*L&\overset{\phi }{\simeq} & f^*p_*\pi^*L&=V\otimes_k{\mathcal O}_{Y'}\\
 & _{\text{Prop 5.1}}\parallel&&\downarrow&&\parallel&\\
 & p_*\bar{f}^*\pi^*L& \overset{\tau}{\simeq}& p_*\bar{f}^*\tau^*\pi^*L&\overset{\phi }{\simeq}& p_*\bar{f}^*\pi^*L\\
  &\parallel&&\parallel&&\parallel&\\
  &p_*\pi^*L& \overset{\tau}{\simeq}& p_*\bar{\tau}^*\bar{f}^*\pi^*L&\overset{\phi }{\simeq}& p_*\pi^*L\\
  &\parallel&&\parallel&&\parallel&\\
 V\otimes_k{\mathcal O}_{Y'}= & p_*\pi^*L& \overset{\bar{\tau}}{\simeq}& p_*\bar{\tau}^*\pi^*L&\overset{\bar{\phi} }{\simeq}& p_*\pi^*L&=V\otimes_k{\mathcal O}_{Y'}\end{matrix}$$

b) Squares up and down of the first diagram are commutative:

\medskip

$\begin{matrix}\tau_1^*\pi^* L &\overset{\tau_2}{\to}& \tau_{2*}\tau_2^*\tau_1^*\pi^* L  \\
\downarrow \phi_1& &\downarrow \phi_1\\
\pi^* L &\overset{\tau_2}{\to}& \tau_{2*}\tau_2^*\pi^* L \\
\downarrow \tau_2& &\downarrow \phi_2\\
\tau_{2*}\tau_2^*\pi^* L &\overset{\phi_2}{\to}& \tau_{2*}\pi^* L \end{matrix}$ \ and taking $\tau_{1 *}$ \  $\begin{matrix}\tau_{1 *}\tau_1^*\pi^* L &\overset{\tau_2}{\to}& \tau_{1 *}\tau_{2*}\tau_2^*\tau_1^*\pi^* L  \\
\downarrow \phi_1& &\downarrow \phi_1\\
\tau_{1 *}\pi^* L & & \tau_{1 *}\tau_{2*}\tau_2^*\pi^* L \\
\downarrow \tau_2& &\downarrow \phi_2\\
\tau_{1 *}\tau_{2*}\tau_2^*\pi^* L &\overset{\phi_2}{\to}& \tau_{1 *}\tau_{2*}\pi^* L \end{matrix}$

\bigskip

Composing  with \  $\pi^* L \overset{\tau_1}{\to} \tau_{1 *}\tau_1^*\pi^* L $ \ and taking $p_*$ one concludes:

$$\begin{matrix}p_*\pi^* L &\overset{\tau_1}{\to} &p_*\tau_1^*\pi^* L &\overset{\tau_2}{\to}& p_*\tau_2^*\tau_1^*\pi^* L  \\
&&\downarrow \phi_1& &\downarrow \phi_1\\
&&p_*\pi^* L & & p_*\tau_2^*\pi^* L \\
&&\downarrow \tau_2& &\downarrow \phi_2\\
&&p_*\tau_2^*\pi^* L &\overset{\phi_2}{\to}& p_*\pi^* L \end{matrix}$$

\medskip

$(\tau_2\circ\tau_1)_{\bar{\phi}}=\phi_2\circ\phi_1\circ\tau_2\circ\tau_1\overset{\text{diagram}}{==}
\phi_2\circ\tau_2\circ\phi_1\circ\tau_1=\tau_2{_{\phi_2}}\circ\tau_1{_{\phi_1}}$.

 \end{pf}
\begin{cor} Let $Y$ be a smooth variety. If \  $\tau\in Aut_Y(U\times Y)$ \  and

$\phi_1,\phi_2:\tau^*L\simeq \pi^*L$ are isomorphisms, then \   $\tau_{\phi_1}=\tau_{\phi_2}\circ h_\lambda$ \  where \ $h_\lambda$ \  is the morphism  multiplying  by \  $\lambda\in {\mathcal O}_Y(Y)^\times$.
\end{cor}
\begin{pf}
 $\phi_{1}\circ\phi_{2}^{-1}$ is an automorphism  of the invertible sheaf $\pi^*L$. Therefore   \  $\phi_{1}\circ\phi_{2}^{-1}=h_{\lambda}$ \ (multiplying by $\lambda$) \  where \  $\lambda\in {\mathcal O}_{U\times  Y}(U\times  Y)^\times={\mathcal O}_Y(Y)^\times$ ( By proposition 2.4). Applying the part  b) of the theorem, $\tau_{\phi_1}=Id_{\phi_1\circ\phi_2^{-1}}\circ \tau_{\phi_2}=h_\lambda \circ \tau_{\phi_2}$

\end{pf}

Let us show what is this morphism when $Y=Spec\,k$ is a  rational point: Let \  $\tau\in Aut_kU$ \  an automorphism  and  \  $\phi: \tau^*L \simeq L$ \  an isomorphism.

We will denote the same way the constant sheaf \   $k(X)=k(U)$.

 If \  $L ={\mathcal O}_U(H )$,  $\tau^*L ={\mathcal O}_U(\tau^{-1}(H ))$, the sections of these sheafs are elements of the fraction field $k(U)$ \  and we have the following commutative diagrams:

$$\begin{matrix}\begin{matrix}{\mathcal O}_U&\overset{\tau}{\to}&\tau_*{\mathcal O}_U\\
 \downarrow& &\downarrow\\
 k(U)&\overset{\tau}{\to}&\tau_*k(U)=k(U)\end{matrix}&\hskip 1cm&\begin{matrix}{\mathcal O}_U(H )&\overset{\tau}{\to}&\tau_*{\mathcal O}_U(\tau^{-1}H )\\ \downarrow& &\downarrow\\ k(U)&\overset{\tau}{\to}&k(U)\end{matrix}\end{matrix}$$

 The isomorphism \  $\phi:\tau^*L\simeq L$ \  induces an isomorphism over $k(U)$ which is multiplying by some  $0\neq f\in k(U)$ \  obtaining the commutative diagrams:
$$\begin{matrix}L&\overset{\tau}{\to}& \tau_*\tau^*L&\overset{\phi}{\simeq}& \tau_*L \\
                    \downarrow& &             \downarrow&                                   &\downarrow\\
k(U)&\overset{\tau}{\to}& k(U)&\overset{\cdot f}{\simeq}&k(U) \end{matrix}$$

Taking global  section one has that is
 $a\in V$ is considered  as an  element of  $ k(U)$, then \   $\tau_\phi(a)=\tau(a)\cdot f$.

Now we have the invertible sheaves \  $L_1={\mathcal O}_U(H_1),\dots, L_r={\mathcal O}_U(H_r)$ \  of $U$. Let \  $\tau: U\times Y\to U\times Y$ \ be  an automorphism  over  $Y$ \ and for each $i$,\  $\phi_i:\tau^*\pi^*L_i\simeq \pi^*L_i$ \  an isomorphism. For each  $n=(n_1,\dots,n_r)\in {\mathbb Z}^r$, let \ $L^n=L_1^{n_1}\otimes\dots \otimes L_r^{n_r}=$

$={\mathcal O}_U(n_1H_1+\dots+n_rH_r)$. These isomorphisms produce  an isomorphism \  $\phi_n: \tau^*\pi^*L^n\simeq \pi^*L^n$  \  and so we have a ${\mathbb Z}^r$-graded algebra isomorphism  \  $\phi:\tau^*\pi^*{\mathcal A}\to \pi^*{\mathcal A}$. This pair \  $(\tau,\phi)$ \  produces an element \
$\tau_\phi\in Aut_{{\mathcal O}_Y-gr}( A\otimes_k{\mathcal O}_Y)$ \  as we have already seen in this section. By theorem 5.2 one has:

\begin{thm} a) The correspondence  \  $(\tau,\phi)\to \tau_\phi$ \  is functorial. That is  if $(\tau,\phi)$ is as above and  \   $f:Y'\to Y$ is a morphism, then \  $f^*(\tau_\phi)=(f^*\tau)_{f^*\phi}$.

b) The previously correspondence is a group morphism :

   Let \  $\tau_1,\tau_2\in Aut_Y(U\times Y)$ \  and  \  $\phi_1:\tau_1^{*}\pi^*{\mathcal A}\to \pi^*{\mathcal A} \ , \ \phi_2:\tau_2^{*}\pi^*{\mathcal A}\to \pi^*{\mathcal A} $ \  be ${\mathbb Z}^r$-graded isomorphisms.
 If \  $\bar{\phi}$ \  is the composition of the isomorphisms \  $\tau_2^*\tau_1^{*}\pi^*{\mathcal A}\overset{\phi_1}{\simeq} \tau_2^*\pi^*{\mathcal A}\overset{\phi_2}{\simeq} \pi^*{\mathcal A}$, then
 \  $(\tau_2\circ\tau_1)_{\bar{\phi}}=\tau_{2\phi_2}\circ \tau_{1\phi_1}$ .

\end{thm}

By corollary 5.3 we deduce:

\begin{cor} Let  $Y$ be a smooth variety. If  \  $\tau=\tau_1=\tau_2$, then  \  $\tau_{\phi_1}=\tau_{\phi_2}\circ h_\lambda$ \  where \ $\lambda\in {\mathcal O}_Y(Y)^{\times r}$. Therefore, $[\tau_{\phi_1}]=[\tau_{\phi_2}]$ \  as element of \  $(Aut_gA/T_{{\mathbb Z}^r})(Y)$.
\end{cor}

 Let \  $Z=Aut_k^0X$ \  the connected component of $Aut_kX$ in the rational point $e=Id$.
 If \ $\tau: X\times Z\to X\times Z$ \  is the universal automorphism, then it induces an  automorphism over  nonsingular points and so an automorphism  \  $\tau:U\times Z\to U\times Z$ \  such that  \  $\tau_{|U\times \{e\}}=Id$.

By proposition 2.2 \  $\tau^*\pi^*L_i\simeq \pi^*L\otimes p^* L_i'$. Restricting to \  $U\times \{e\}$ \  one has that  \  $\pi^* L_i\simeq \pi^* L$ \  and  \  so \  $\tau^*\pi^*L_i\simeq \pi^*L_i\otimes p^* L_i'$.

There exists an  affine open covering \  $Z=\underset{l}{\cup} U_l$ \  in such a way that for all $i, l$, $L'_{i\ |U_l}\simeq {\mathcal O}_{U_l}$. Let  \  $R=\underset{l}{\coprod} U_l\to Z$ \  this covering. $f:U\times R\to U\times Z$ is a open covering and we denote  $f^*\tau=\tau$ since $f^*\tau$ is only to  restrict  $\tau$ to each open set of the covering.  On  \  $U\times R$, we have an isomorphism \  $\phi_i:\tau^*\pi^*L_i\simeq \pi^*L_i$ on $U\times R$. Therefore we have a ${\mathbb Z}^r$-graded isomorphism \ $\phi:\tau^*\pi^*{\mathcal A}\simeq \pi^*{\mathcal A}$ \  on $U\times R$ \   and so a graded automorphism,  $\tau_{\phi}: A\otimes_k{\mathcal O}_{U\times R }\to A\otimes_k{\mathcal O}_{U\times R }$. Let $R\times_ZR=\underset{l,l'}{\coprod} \, U_l\cap U_{l'}$ \   the intersection covering and
 \  $p_1,p_2:U\times R\times_Z R\rightrightarrows U\times R$ \   the canonical projections. The isomorphism $\phi$ gives isomorphisms \   $\phi_{1}:p_1^*\tau^*\pi^*{\mathcal A}\overset{\phi}{\simeq} p_1^*\pi^*{\mathcal A}=\pi^*{\mathcal A}$ \  and \  $\phi_{2}: p_2^*\tau^*\pi^*{\mathcal A}\overset{\phi}{\simeq} p_2^*\pi^*{\mathcal A}=\pi^*{\mathcal A}$. As $\pi\circ\tau\circ p_1=\pi\circ\tau\circ p_2$ \  we obtain that \ $\phi_{1}\circ \phi_{2}^{-1}$ \ is an automorphism of \  $\pi^*{\mathcal A}$ (on $U\times R\times_Z R$) \ and so  \
 $\phi_{1}\circ \phi_{2}^{-1}=h_{\lambda}$ \  where \  $\lambda\in {\mathcal O}_{R\times_Z R}(R\times_Z R)^{\times r}$. By the corollary 5.3, $\tau_{\phi_1}=\tau_{\phi_2}\circ h_{\lambda}$. Therefore \  $[\tau_{\phi_1}]=[\tau_{\phi_2}]$ \  in  \  $Aut (A\otimes_k{\mathcal O}_{R\times_Z R})/{\mathcal O}_{R\times_Z R}(R\times_Z R)^{\times r}$.

That is to say, we have an element \  $\varphi=[\tau_\phi]$ \  of  $Aut^0_gA/T_{{\mathbb Z}^r} \,(Z)$.

The morphism defined for any $k$-scheme $Y$ is: If $g:Y\to Z$ is a morphism; we have  an open cover of $Y$, $R_Y=R\times_ZY\to Y$, and a morphism  $g':R_Y\to R$. Therefore \  $g'{^*}\tau_\phi\in Aut_{{\mathcal O}_{R_Y}-gr}(A\otimes_k{\mathcal O}_{R_Y})$ \  which coincides over intersections (mod $T_{{\mathbb Z}^r}$) and so  produces an element of \  $Aut_g^0A/T_{{\mathbb Z}^r} \,(Y)$.

\begin{thm} Let  \  $\varphi: Z\to Aut^0_gA/T_{{\mathbb Z}^r} $ \  be the above constructed morphism.

a) $\varphi$ \ does not  depend on  the open covering  \   $R\to Z$ \  or the isomorphism
\  $\phi:\tau^*\pi^*{\mathcal A}\simeq \pi^*{\mathcal A}$ \  on $U\times R$.

b) $\varphi $ \  is a group morphism.

c)  $\varphi: Z=Aut_k^0X\to Aut_gA/T_{{\mathbb Z}^r} $ \  is injective.

d) The induced morphism between the tangent spaces \  $\varphi_{*}: T_eAut_k^0X\to T_eAut_gA/T_{{\mathbb Z}^r} $ \  is an isomorphism.

e) $Aut^0_kX$ is isomorphic to \  $Aut^0_gA_C \, / \, T_{N}$.
\end{thm}

 \begin{pf} a) Let \  $f:R\to Z$ \  and \  $f':R'\to Z$ \ be two open coverings of $Z$ \  with isomorphisms \   $\phi:\tau^*\pi^*{\mathcal A}\to \pi^*{\mathcal A}$ \  on $U\times R$ \  and \  $\phi':\tau^*\pi^*{\mathcal A}\to \pi^*{\mathcal A}$ \   on  $U\times R'$. Let \  $\bar{R}=R\times_ZR'\to Z$ \ the intersection covering and \  $p_1:\bar{R}\to R$, $p_2:\bar{R}\to R'$ \  the canonical projections. We have that  \   $p_1^*\phi$ \ and \  $p_2^*\phi'$ \  are isomorphisms from  $\tau^*\pi^*{\mathcal A}$ into  $\pi^*{\mathcal A}$ \  on   $U\times\bar{R}$. Therefore  \ $p_1^*\phi=p_2^*\phi'\circ h_\lambda$ \  where \  $\lambda\in {\mathcal O}_{\bar{R}}(\bar{R})^{\times r}$.

  So \  $\tau_{p_1^*\phi}=\tau_{p_2^*\phi'}\circ h_\lambda$ \  and one has that  \  $p_1^*\tau_\phi\equiv p_2^*\tau_{\phi'}$  \  mod $T_{{\mathbb Z}^r}(\bar{R})$. Therefore  \  $[\tau_\phi]=[\tau_{\phi'}]$ \ on $Aut_gA/T_{{\mathbb Z}^r}\,(Z)$ \

b) If  \  $\tau_1,\tau_2\in Aut_Y(X\times Y)$, then  there are  open coverings  $R_1$ and $R_2$ of $Y$ where \  $\tau_1^*\pi^*{\mathcal A}$ and $\tau_2^*\pi^*{\mathcal A}$ \ are isomorphic  to  \  $\pi^*{\mathcal A}$. Taking the intersection covering \   $R=R_1\times_YR_2$, one has isomorphisms on $U\times R$ from  \  $\tau_1^*\pi^*{\mathcal A}$ and $\tau_2^*\pi^*{\mathcal A}$ \  into  $\pi^*{\mathcal A}$ \  which produce
automorphisms of \  $ A\otimes_k{\mathcal O}_R$. By  theorem 5.4 b) one can conclude.

c) We have to prove that \  $Ker\,\varphi=0$. We will prove it on the rational points. Let $\tau$ be an automorphism of $X$ and $\phi:\tau^*{\mathcal A}\simeq {\mathcal A}$ an isomorphism. We assume that $\tau_\phi=h_\lambda$. The ${\mathbb Z}^r$-graded automorphism  $\tau_\phi:A\to A$ induces an automorphism on $k(X)$ (Proposition 4.3 a)). It is enough to prove this automorphism is precisely  the automorphism that $\tau$ induces on $k(X)$.
 Let us show it: There is an injective morphism \  ${\mathcal O}_U(H_i)=L_i\to k(X)$ \ which
 is isomorphism at generic point. Therefore there is an injective ${\mathbb Z}^r$-graded  morphism \  $ {\mathcal A}\to k(X)[y_1,\dots,y_r]$ \  given by  \  $L^n\to k(X)\cdot y^n$ \  which is isomorphism at generic point.

We know that the isomorphism  \  $\phi_i:\tau^*L_i\simeq L_i$ \  induces an isomorphism on the fraction field
which is multiplying by some $0\neq f_i\in k(X)$. So the isomorphism \  $\phi_n:\tau^*L^n\simeq L^n$, for each $n=(n_1,\dots,n_r)$ induces an isomorphism  on the fraction field which is multiplying by \   $ f^n=f_1^{n_1}\cdot{\dots}\cdot f_r^{n_r}\in k(X)$.

Therefore we have the commutative diagrams:

$$\begin{matrix}{\mathcal A}&\overset{\tau}{\to}& \tau_*\tau^*{\mathcal A}&\overset{\phi}{\simeq}& \tau_*{\mathcal A} \\
                    \downarrow& &             \downarrow&                                   &\downarrow\\
k(X)[y_1,\dots,y_r]&\overset{\tau}{\to}& k(X)[y_1,\dots,y_r]&\overset{h_f}{\simeq}&k(X)[y_1,\dots,y_r] \end{matrix}$$

Taking global sections, one has that if \ $a\in A_n$ is considered as an element of  $ k(X)$, then
$\tau_\phi(a)=\tau(a)\cdot f^n$.
By the proposition 4.3, if $\frac{p}{q}\in k(X)$ where $p,q\in A$ are the same degree $n$, then \  $\tau_{\phi}(\frac{p}{q})=\frac{\tau_{\phi}(p)}{\tau_{\phi}(q)}=\frac{\tau(p)\cdot f^n}{\tau(q)\cdot f^n}=\tau(\frac{p}{q})$.

If $\tau_\phi=h_\lambda$ for some $\lambda\in k^{\times \, r}$, then $\tau=Id$ on $k(X)$ and so $\tau=Id$ on $X$.

d) By c) one just  has  to prove that the tangent spaces are the same dimension. By the corollary 4.7 we only have to prove that $Aut_k^0X$ and $Aut^0_gA_C \, / T_N$ have the same number of roots and this is true by the remark 4.8.

e) $\varphi$ is an isomorphism by proposition 1.1 and one concludes by corollary 4.7.
\end{pf}

\medskip\section{Automorphisms of Cox ring}

We denote $R_u$ to the unipotent radical of $Aut_g^0A_C$ and $G_R=Aut_g^0A_C/R_u$ to the reductive group associated to $Aut^0_gA_C$.

 $\Delta_1=\{v_1,\dots ,v_r\}$ and $A_C$ as the section 4. We can also think  $\Delta_1$ as the set of the rays of the toric variety or the set of the $T_M$-invariant hypersurfaces of the toric variety or the set of variables of $A_C$.

We know, by the section 4.1, that the roots of $Aut^0 _gA_C$ are
\  $r_{ij}\in {\mathbb Z}^r$ \  holding that $\pi(r_{ij})=0$ , $v_i(r_{ij})=-1$ \ and \  $v_j(r_{ij})\geq 0$ for all $j\neq i$.

 \begin{dfn} We will say that the root $r$ is a root associated to $v_i$ if $v_i(r)=-1$. A root $r$ is said to be semisimple if $-r$ is a root. We will say that the semisimple root $r_{ij}$ is  associated to $(v_i,v_j)$ if $r$ is associated to $v_i$ and $-r$ is associated to $v_j$. That is,  $v_i(r)=-1$, $v_j(r)=1$ and $v_k(r)=0$ for all $k\neq i,j$.
 \end{dfn}

 \begin{prop} If \   $\tau:A_C=k[y_1,\dots,y_r]\to A_C=k[y_1,\dots,y_r]$  \  is a  $N$-graded morphism, then \   $\tau(y_j)=y_j(\underset{i}{\sum} \lambda_i y^{r_{ij}})$  where $r_{ij}$ are roots associated to $v_j$.
 \end{prop}
 \begin{pf} $\tau(y_j)$ has the same degree as $y_j$. But
 $y^n$ has the same degree as $y_j=y^{e_j}$ if and only if \  $n=e_j+v(\alpha)\geq 0$. Therefore $v(\alpha)=r_{ij}$ is a root of $Aut_gA_C$ associated to $v_j$ and $y^n=y_jy^{r_{ij}}$.
 \end{pf}

 \begin{prop} Let $r_i,r_j$ be  roots associated to $v_i$ and $v_j$ respectively. We suppose that \  $v_i\neq v_j$ \  and \  $r_i+r_j\neq 0$.

  If \  $v_j(r_i)>0$, then  $v_i(r_j)=0$ \  and \ $r_i+s\cdot r_j$ is a root associated to $v_i$ for each \ $0\leq s\leq v_j(r_i)$. Besides if $r_i$ is not semisimple so  is $r_i+sr_j$.

\end{prop}
\begin{pf}   $v_k(r_i+sr_j)=v_k(r_i)+sv_k(r_j)\geq 0$ for all $v_k\neq v_i,v_j$.

$v_j(r_i+sr_j)=v_j(r_i)-s\geq 0$ \  and \  $v_i(r_i+sr_j)=-1+sv_i(r_j)\geq -1$. If \  $v_i(r_i+r_j)\geq 0$, then $r_i+r_j=0$. Therefore \   $-1=v_i(r_i+r_j)=-1+v_i(r_j)$ \  and we can conclude that  $v_i(r_j)=0$ \  and \ $v_i(r_i+sr_j)=-1$.

We prove that \  $r_i+sr_j$ \  is not semisimple for $s>0$: There exists \  $v_k\in \Delta_1$  \  such that \  $v_k(r_j)>0$. We know that $v_k\neq v_i,v_j$. Therefore $v_k(r_i+sr_j)=v_k(r_i)+sv_k(r_j)>0$. If $r_i+sr_j$ is a semisimple root, it is a root associated to $(v_i,v_k)$. As \  $v_k(r_i)+sv_k(r_j)=1$, we obtain that \  $v_k(r_i)=0$ \  and \ $s=1$. Besides \ $v'(r_i+sr_j)=v'(r_i)+sv'(r_j)=0$ \  for all \  $v'\neq v_i,v_k$. Hence \  $v_j(r_i)=1$,  and $v'(r_i)=0$ for $v'\neq v_j,v_i$. We obtain that $v'(r_i+r_j)\geq 0$ for all $v'\in \Delta_1$ and so $r_i+r_j=0$.

\end{pf}

The   semisimple roots give an equivalence relation on $\Delta_1=\{v_1,\dots,v_r\}$.

\begin{dfn} We say that $v_i\equiv v_j$ if \  $v_i=v_j$ or there is a semisimple root $r_{ij}$ associated to $(v_i,v_j)$.
\end{dfn}

\begin{prop} $v_i\equiv v_j  \  \Leftrightarrow \ y_j=y_iy^{r_{ij}} $ where $r_{ij}=0$ or $r_{ij}$ is a semisimple root associated to $(v_i,v_j) \   \Leftrightarrow \ \text{degree of } y_i=\text{degree of }y_j \   \Leftrightarrow \ H_i\equiv H_j$ modulo  linear equivalence.
\end{prop}
\begin{pf} $r_{ij}$ is a root associated to $(v_i,v_j)  \   \Leftrightarrow \ y_j=y_iy^{r_{ij}} \text{ and }  r_{ij}=v(\alpha_{ij})  \   \Leftrightarrow \ \text{degree of } y_i=\text{degree of }y_j  \   \Leftrightarrow \ H_j-H_i=div \ x^{\alpha_{ij}}$.
\end{pf}

\begin{subsection}{Reductive subgroup of \   $Aut^0_gA_C $}

 \

 \medskip

 All roots of a reductive group  are semisimple ([Hu] or [Sa]).

Let $F\subseteq \Delta_1$ an equivalence class and $A_F=k[y_{i}]_{v_i\in F}$ and we denote \  $Gl_F$ \ as the linear group of order \  $|F|$ \  acting in a natural way on \
$A_F$ (Section 2. d) )   and  by the identity on the rest variables, it acts on \  $A_C$. Therefore \  $Gl_F$ \  is a subgroup of \  $Aut_gA_C$.
 Besides, for $F\neq F'$ the elements of $Gl_F$ commute with the elements of $Gl_{F'}$ (inside  \  $Aut_gA_C$).
Therefore we have that \  $Gl_{F_1}\times\dots\times Gl_{F_k}$ is a subgroup of  $Aut_g^0A_C$.

\begin{thm} $Gl_{F_1}\times\dots\times Gl_{F_k}\simeq G_R$.
\end{thm}
\begin{pf}The maximal torus $T_{{\mathbb Z}^r}$ is contained in $G_R$ because it coincides with the group of diagonal matrices in each $Gl_{F_i}$.
To show that  \  $G_R$ \  coincides with  the reductive group of  \  $Aut^0_gA_C$ \  it is sufficient to prove that they  have the same  tangent space  at  $e$ and for this it is enough to prove that the  roots  of $G_R$ are the semisimple roots of $Aut_g^0A_C$.
 If $r_{ij}$ is a semisimple root   associated to $(v_i,v_j)$, then $y^{r_{ij}}=\frac{y_j}{y_i}$. By the  corollary 3.7,  $T_e(Gl_{F_1}\times\dots\times Gl_{F_k})$ inside \ $T_e(Aut_gA_C)=Der_{N-gr}(A_C)$ \ is
generated by the derivations \   $y_j\frac{\partial}{\partial y_i}=y^{r_{ij}}D_{v_i}$ and we can conclude.
\end{pf}
\end{subsection}
\begin{subsection}{Unipotent radical of $Aut^0_gA_C$}

\

\medskip

By the above section $T_eR_u$ is generated by $y^{r_i}D_{v_i}$ where $r_i$ are the no semisimple roots.
 We denote \  $W=\Delta_1/\equiv \ =\{F_1,\dots,F_k\}$ \  the quotient of $\Delta_1$ modulo $\equiv$.

\begin{dfn} We say  that \  $v_i<v_j$ \  if there is  a no semisimple root $r_i$ associated to $v_i$ such that  $v_j(r_i)>0$.
\end{dfn}

\begin{prop} a)  $<$ \  is a partial order on $\Delta_1$.

  b) Let $F$ and $F'$ be two equivalence class. If  $v_i<v_j$ for some $v_i\in F$ and $v_j\in F'$, then $v_k<v_s$ for all $v_k\in F$ and $v_s\in F'$.

c) The relation on $W$ given by $F<F'$ if there is $v_i\in F$ and $v_j\in F'$ such that $v_i<v_j$ \  is a partial order.

\end{prop}

\begin{pf}
a) If $v_i<v_j$ there is a no semisimple root $r_i$ associated to $v_i$ such that $v_j(r_i)>0$. If $r_j$ is a root associated to $v_j$, by the  proposition 6.3,  $v_i(r_j)=0$ and so it is not possible that $v_j<v_i$.

    Let us show that if $v<v'$ and $v'<v''$, then $v<v''$. If  $v<v'$ and $v'<v''$, there are no semisimple roots $r$ and $r'$ associated respectively to $v$ and $v'$ such that $v'(r)>0$ and $v''(r')>0$. By the proposition 6.3,  $r+r'$ is a no semisimple root associated to $v$  \  and \  $v''(r+r')=v''(r)+v''(r')>0$. Hence \  $v<v''$.

b) There is a no semisimple root $r_i$ associated to $v_i$ such that $v_j(r_i)>0$.

 If $v_k\equiv v_i$, there is a semisimple root $r_{ki}$ associated to $(v_k,v_i)$. As \  $v_i(r_{ki})=1>0$, by the proposition 6.3 we have that \   $v_k(r_i)=0$.  We obtain that  $r_i+r_{ki}$ is a  root associated to $v_k$. If it were a semisimple root associated to $(v_k,v_s)$, then \ $r_i=(r_i+r_{ki})-r_{ki}$  would be a semisimple root associated $(v_i,v_s)$. Therefore $r_i+r_{ki}$ is a no semisimple root associated to $v_k$ and $v_j(r_i+r_{ki})>0$  \  and so \  $v_k<v_j$.

 If $v_j\equiv v_k$, there is a semisimple root $r_{jk}$ associated to $(v_j,v_k)$.  By the proposition 6.3,   $r_i+r_{jk}$ is a no semisimple root associated to $v_i$ and  $v_k(r_i+r_{jk})=v_k(r_i)+1> 0$.  Therefore $v_i<v_k$.

 c) If $F<F'$, there is $v_i\in F$ and $v_j\in F'$ such that $v_i<v_j$. By b) $v<v'$ for all  $v\in F$ and $v'\in F'$. By a) $F'\nless F$.

     If $F<F'$ and $F'<F''$,  $v<v'$ and $v'<v''$ for all  $v\in F$, $v'\in F'$ and $v''\in F''$.  Therefore $v<v''$ and so $F<F''$.
\end{pf}

\begin{subsection}{ Automorphisms given by derivations}

 Let  \ $r_1,\dots,r_l$ \  be all no semisimple roots associated to $v_i$. That is $v_i(r_j)=-1$ and $v'(r_j)\geq 0$ for all $v'\neq v_i$.
 We denote  \  $t=(t_1,\dots,t_l)$ \  and \  $tD=\overset{l}{\underset{j=1}{\sum}}\, t_j y^{r_j}D_{v_i}$.

 $\tau_{tD}:A_C[t_1,\dots,t_l]\to A_C[t_1,\dots,t_l]$  \  is the $N$-graded morphism and $k[t_1,\dots,t_l]$-morphism given by \
 $\tau_{tD}(y^n)=y^n(1+t_1\cdot y^{r_1}+\dots+t_l\cdot y^{r_l})^{v_i(n)}$ \  for each \  $n\in {\mathbb N}^r$.

 \begin{rem} $y^n(1+t_1\cdot y^{r_1}+\dots+t_l\cdot y^{r_l})^{v_i(n)}=$

 $=\underset{k_1+\dots+k_l\leq v_i(n)}{\sum} \,\lambda_{k_1,\dots,k_l}t_1^{k_1}\cdot t_2^{k_2}\cdot{\dots}t_l^{k_l}\cdot y^{n+k_1r_1+\dots+k_lr_l}\in k[{\mathbb N}^r][t_1,\dots,t_l]$ \  because

 $v_i(n+k_1r_1+\dots+k_lr_l)=v_i(n)-k_1-\dots-k_l\geq 0$

 \end{rem}

It is a ring morphism because  $\tau_{tD}(y^n\cdot y^m)=\tau_{tD}(y^n)\cdot \tau_{tD}(y^m)$ and it is an automorphism because the following proposition proves that $\tau_{tD}^{-1}=\tau_{-tD}$.

By the proposition 4.3 b) and theorem 5.6 e), these automorphisms produce  automorphisms in $X$ and $k(M)$. These are: If \   $tD=\overset{l}{\underset{j=1}{\sum}}\, t_j x^{\alpha_j}D_{v_i}\in Der({\mathcal O}_X)$, $\tau_{tD}(x^\alpha)=x^\alpha(1+t_1\cdot x^{\alpha_1}+\dots+t_l\cdot x^{\alpha_l})^{v_i(\alpha)}$ \  for each \  $\alpha\in M$.

 \begin{prop}a) For \   $t=0=(0,\dots,0)$, $\tau_0=Id$ \ and

  $\tau_{(t+t')D}=\tau_{tD}\circ \tau_{t'D}$ \ being \ $t=(t_1,\dots,t_l)$ and $t'=(t'_1,\dots,t'_l)$.

b) The morphism \  $f_{v_i}:{\mathbb A}^1\times\dots\times {\mathbb A}^1=Spec\,k[t_1,\dots,t_l]\to Aut_gA_C$ \  given by \  $f(t)= \tau_{tD}$ \  is an injective morphism of groups.

c) The tangent space on identity element of the subgroup $Im\,f_{v_i}$ \ is the $k$-vectorial space  generated by \ $\{y^{r_1}D_{v_i},\dots,y^{r_l}D_{v_i}\}$.

d) $Im\,f_{v_i}$ is an additive subgroup of $R_u$.
 \end{prop}
 \begin{pf} a) We denote  \  $t\cdot y^r=1+t_1\cdot y^{r_1}+\dots+t_l\cdot y^{r_l}$.

 $\tau_{tD}(y^{r_j})=y^{r_j}{(t\cdot y^r)}^{v_i(r_j)}=\dfrac{y^{r_j}}{t\cdot y^r} $.

$\tau_{tD}\circ \tau_{t'D}(y^n)=\tau_{tD}(y^n(1+t'_1\cdot y^{r_1}+\dots+t'_l\cdot y^{r_l})^{v_i(n)})=$

$\tau_{tD}(y^n)(1+t'_1\cdot\tau_{tD}(y^{r_1})+\dots+t'_l\cdot \tau_{tD}(y^{r_l}))^{v_i(n)}=$

$y^n{(t\cdot y^r)}^{v_i(n)}(1+t'_1\cdot \dfrac{y^{r_1}}{t\cdot y^r}+\dots+t'_l\cdot \dfrac{y^{r_l}}{t\cdot y^r})^{v_i(n)}=$

$y^n{(t\cdot y^r)}^{v_i(n)}\dfrac{(t\cdot y^r+t'_1\cdot y^{r_1}+\dots+t'_l\cdot y^{r_l})^{v_i(n)}}{{(t\cdot y^r)}^{v_i(n)}}=y^n(1+(t+t')\cdot y^r)^{v_i(n)}=\tau_{(t+t')D}(y^n)$.

b) a) proves that $f_{v_i}$ is a morphism of groups. It is an injective morphism if \   $Ker\,f_{v_i}=\{0\}$. If \  $\tau_{tD}=Id$, then $\tau_{tD}(y^n)=y^n(1+t_1\cdot y^{r_1}+\dots+t_l\cdot y^{r_l})^{v_i(n)}=y^n$. Taking $n\in {\mathbb N}$ such that $v_i(n)=1$,   $1+t_1\cdot y^{r_1}+\dots+t_l\cdot y^{r_l}=1$. Hence $t_1=t_2=\dots=t_l=0$.

c) Let $tD=ty^r D_{v_i}$ and  let  $f_{v_i}:{\mathbb A}^1=Spec\,k[t]\to Aut_gA_C$ be the morphism induced. It sufficient to prove   that the induced morphism in space tangent at $e$, $f_{v_i\,*}:T_0{\mathbb A}^1\to T_e(Aut_gA_C) $ \  satisfies that \   $f_{v_i\,*}((\frac{\partial}{\partial t})_0)=y^rD_{v_i}$.

 The derivation $(\frac{\partial}{\partial t})_0$ produces the morphism \  $g:Spec\,k[\epsilon]\to {\mathbb A}^1=Spec\,k[t]$ \  defined by   \  \  $\begin{matrix}k[t]&\to& k[\epsilon]\\ p(t)&\to & p(0)+p'(0)\epsilon\end{matrix}$ . In particular $t\to\epsilon$.

$f_{v_i\,*}((\frac{\partial}{\partial t})_0)$ \ is the derivation  given by the  morphism \   $f_{v_i}\circ g:Spec\,k[\epsilon]\to Aut_gA_C$ \  whose corresponding $k[\epsilon]$-automorphism is \   $A_C[t]\otimes_{k[t]}k[\epsilon]\overset{\tau_{tD}\otimes 1}{\to} A_C[t]\otimes_{k[t]}k[\epsilon]$. This automorphism is such that \   $y^n\to y^n(1+y^r\epsilon)^{v_i(n)}\overset{\epsilon^2=0}{=}y^n(1+v_i(n)y^r\epsilon)=y^n+Dy^n \, \epsilon$. That is  the $k[\epsilon]$-automorphism produced by the derivation $D$.

d) $Im\,f_{v_i}\subset R_u$ \  because \ $T_e(Im\,f_{v_i})\subset T_e(R_u)$ (proposition 1.1).
 \end{pf}

 We denote $V_i=Im\, f_{v_i}$. In the case that $v_i$ has not associated any root or it has  only associated semisimple roots, $V_i=0$. We can think $V_i$ as \underline{the  underlying group} \underline{under addition of the $k$-vectorial space generated by derivations} \  $\{y^{r_1}D_{v_i},\dots,y^{r_l}D_{v_i}\}$.

   \begin{prop} Let $D=y^{r_i}D_{v_i}$ and $D'=y^{r_j}D_{v_j}$   be  graded derivations of $A_C$.
   \begin{enumerate}
   \item[a)] If  $v_i(r_j)=0=v_j({r_i})$,  $\tau_{tD}$ and $\tau_{t'D'}$ commute.

   \item[b)] If $v_j({r_i})=0$, $\tau_{tD}\circ \tau_{t'D'}\circ \tau_{tD}^{-1}=\tau_{\bar{t}\, D'}$ \  where \  $\bar{t}\,D'=\overset{v_i(r_j)}{\underset{k=1}{\sum}} \bar{t}_ky^{kr_i+r_j}D_{v_j}$.
   \end{enumerate}

\end{prop}

  \begin{pf} a)
   $\tau_{tD}(y^{r_j})=y^{r_j}$ ; $\tau_{t'D'}(y^{{r_i}})=y^{{r_i}}$ ;

  Therefore \  $\tau_{tD}\circ\tau_{t'D'}(y^n)=\tau_{tD}(y^n(1+t'\cdot y^{r_j})^{v_j(n)})=$

  $y^n(1+t\cdot y^{r_i})^{v_i(n)}(1+t'\cdot y^{r_j})^{v_j(n)}$.

  In the same way \   $\tau_{t'D'}\circ\tau_{tD}(y^n)=y^n(1+t'\cdot y^{r_j})^{v_j(n)}(1+t\cdot y^{r_i})^{v_i(n)}$  \  and one concludes.

b)  $\tau_{t'D'}(y^{{r_i}})=y^{{r_i}}$,  $\tau_{tD}(y^{{r_i}})=\dfrac{y^{{r_i}}}{1+ty^{r_i}}$, $\tau_{t'D'}(y^{{r_i}})=y^{{r_i}}$ \ and \  $ \tau_{tD}^{-1}= \tau_{-tD}$.

 $\tau_{tD}\circ \tau_{t'D'}\circ \tau_{tD}^{-1}(y^n)= \tau_{tD}\circ \tau_{t'D'}(y^n(1-ty^{r_i})^{v_i(n)})=$

 $\tau_{tD}( y^n(1+y^{r_j})^{v_j(n)}(1-ty^{r_i})^{v_i(n)})=$

$=y^n(1+ty^{r_i})^{v_i(n)}(1+t'y^{r_j}(1+ty^{r_i})^{v_i(r_j)})^{v_j(n)}(1-t\frac{y^{{r_i}}}{1+ty^{r_i}})^{v_i(n)}=$

$y^n(1+t'y^{r_j}(1+ty^{r_i})^{v_i(r_j)})^{v_j(n)}=y^n(1+t'y^{r_j}(\overset{v_i(r_j)}{\underset{k=0}{\sum}} \begin{pmatrix}v_i(r_j)\\ k\end{pmatrix}t^k y^{k{r_i}}))^{v_j(n)}$.

Putting \  $\bar{t}_k=t'\begin{pmatrix}v_i(r_j)\\ k\end{pmatrix}t^k$, we can conclude because, by proposition 6.3, $r_j+kr_i$ is a no semisimple root associated to $v_j$.

\end{pf}
\begin{cor}a) If $v_i\nless v_j$ and $v_j\nless v_i$, then every element of $V_i$ commutes with every element of $V_j$.

b) If $v_i\equiv v_j$, then every element of $V_i$ commutes with every element of $V_j$.

\end{cor}
\begin{pf} a) In this case $v_i(r')=0=v_j(r)$ for all roots $r,r'$ associated to $v_i$ and $v_j$ respectively.

b) If $v_i\equiv v_j$, then $v_j\nless v_i$ and $v_i\nless v_j$.
\end{pf}
\end{subsection}

For each equivalence class $F$, we denote \  $V_F=\underset{v_i\in F}{\prod} V_i$ \  which is an additive subgroup of $Aut_gA_C$ whose roots are  all no semisimple roots associated to some $v_i\in F$.

 We call dimension   of $z\in W$ and denote by $d(z)$ the maximum of the  length of the chains of $W$ finishing on $z$: $d(z)= {Max} \, \{l\}_{ z_1<\dots<z_l<z}$.

For each $i\geq 0$, we call  $W_i=\{z\in W \, /\ d(z)=i\}$.

 One has that  if $i\leq j$ and $F\in W_i$ and $F'\in W_j$,  then  $F'\nless F$.

 Therefore by the corollary 6.12 a) if $F,F'\in W_i$, $V_F$ commute with $V_{F'}$ \  and so \
  $E_i=\underset{F\in W_i}{\prod} V_F$ \   is an additive subgroup of $R_u$.

\begin{lem} If \  $j\leq i$ and \   $g\in E_i$, then \  $g\cdot E_j\cdot g^{-1}=E_j$.
\end{lem}
\begin{pf} It is enough to prove that for \ $tD=ty^rD_{v_k}$ ( $v_k\in F\in W_i$) and $t'D'=t'y^{r'}D_{v'}$ ( $v'\in F'\in W_j$), $\tau_{tD}\circ \tau_{t'D'}\circ ({\tau_{tD}})^{-1}\in E_j$. As $v_k\nless v'$,  $v'(r)=0$ \  and we can conclude by the proposition 6.11 b).
\end{pf}
 \begin{thm} Let  $R_u$ be the unipotent radical of \  $Aut_{g}^0A_C$.

  a) For all $i$, $L_i=E_1\cdot  E_2\cdot \dots{\cdot}  E_i$ \ is a subgroup of $R_u$ .

 b) $L_l= E_1\cdot  E_2\cdot \dots{\cdot}  E_l$ \   is isomorphic to $R_u$.

 c) $L_i$ is  a normal subgroup of $R_u$ \ and \  $L_i/L_{i-1}$ \  is an additive group.

\end{thm}

\begin{pf} a) We prove it by induction over $i$. For $i=1$, $L_1=E_1$ and there is nothing to say.
 Suppose it is true for $i-1$ and we prove it for $i$.

 $L_{i-1}$ is a subgroup by hypothesis of induction. $L_i=L_{i-1}\cdot E_i$ \  is a subgroup because, by above lemma, $g\cdot L_{i-1}\cdot g^{-1}=L_{i-1} $ \  for every $g\in E_i$.

 b) The tangent space of both groups at $e$   coincides.

 c) $L_i$ is normal in $L_l$ because by above lemma \  $g\cdot L_i\cdot g^{-1}=L_i$ \  for each $i\leq j$ and $g\in E_j$. Besides \  $E_i\simeq L_i/L_{i-1}$ \  and so it is additive group.

 \end{pf}

\begin{cor} $R_u$ is a semidirect product of the  additive groups.
\end{cor}
\end{subsection}

\begin{subsection}{Action of $G_R$ on $R_u$}

\

\medskip

We have shown that each equivalence class $F$ defines an additive group $V_F$ and a linear group $Gl_F$; both subgroups of $Aut_g^0A_C$.

\begin{rem} Both groups, $Gl_F$ and $V_F$, are the identity  over the variables of $A_C$ which are not in $F$.
\end{rem}

Let \  $F=\{v_1,\dots,v_l\}$ \  and \  $E_F=<y_1,\dots,y_l>_{y_i\in F}$ \ the $k$-vectorial space generated by the variables in $F$. We have that \  $Gl_F=Aut_kE_F$.  It is known that \  $E_F^*$ (dual of $E_F$) \  and  \ $S^nE_F$ ($n$-th symmetric power of $E_F$ or homogeneous polynomials of degree $n$ in the variables $y_1,\dots,y_l$) \  are irreducible representations of $Gl_F$.

\begin{thm}
 a) $Gl_F$ acts on $V_F$ and as \   $GL_F$-module is \  $V_F=\oplus E_F^*$.

 b) If $F'<F$, $Gl_F$ acts on $V_k$ for each $v_k\in F'$ \  and \  $V_k={\underset{n_i}{\oplus}} S^{n_i}(E_F)$ as $Gl_F$-module.

c)  If $F'\nleq F$,  $Gl_F$ acts on $V_{F'}$ trivially.

In particular $Gl_F$ acts on $V_{F'}$ for all $F'\in W$.
\end{thm}
\begin{pf} Let \  $F=\{v_1,\dots,v_l\}$,  $E_F=<y_1,\dots,y_l>$ \  and \  $r_{ij}$ \  is the semisimple root associated to $(v_i,v_j)$. That is \ $y^{r_{ij}}=\frac{y_j}{y_i}$.

 a) If $r$ is a no semisimple root associated to $v_j\in F$, then by the proposition 6.3, $r+r_{ij}$   \   is a no semisimple root associated to $v_i$.
 Hence \  dim\,$V_j=$ dim\,$ V_i$. If   the number of no semisimple roots associated to any $v_i\in F$ is $h$, then \  dim $V_F=l\cdot h$.

 Let \  $C=(c_{ij})\in Gl_F$ \  an invertible square matrix of order $l$. The automorphism defined by $C$  will be denoted by $\tau_C$. We have that \      $\tau_C(y_j)=\underset{v_i\in F}{\sum} c_{ij}y_i$.

Let $r $  be  a no semisimple root associated to $v_k\in F$ and \   $tD=t y^r D_{v_k}$.    One has that $v_i(r)=0$ for each $v_i\in F$ $(i\neq k)$ \ and \ $v_k(r)=-1$. Therefore  in $y_ky^r$ does not occur any  variable of  $F$ and so $\tau_C(y_ky^r)=y_ky^r$.

For each $y_i$ such that $v_i\in F$, $\tau_{tD}(y_i)=y_i$ ( $i\neq k$ ) and $\tau_{tD}(y_k)=y_k+ty^{r}y_k$.

Let $C^{-1}=(b_{ij})$ . $\tau_C\circ\tau_{tD}\circ\tau_{C^{-1}}(y_j)=\tau_C\circ\tau_{tD}(\underset{v_i\in F}{\sum}
b_{ij}y_i)=$

$\tau_C(\underset{v_i\in F,i\neq k}{\sum} b_{ij}y_i)+b_{kj}\tau_C\tau_{tD}(y_k)=\tau_C(\underset{v_i\in F,i\neq k}{\sum} b_{ij}y_i)+b_{kj}\tau_C(y_k+ty_k\cdot y^r)=$

$\tau_C(\underset{v_i\in F}{\sum} b_{ij}y_i)+tb_{kj}\tau_C(y_k\cdot y^r)=\tau_C\circ\tau_{C^{-1}}(y_j)+tb_{kj}y_ky^r=y_j+tb_{kj}y_jy^{r_{jk}}y^r=
y_j(1+tb_{kj}y^{r_{jk}+r})=\tau_{t\bar{D}}(y_j)$ \  where \  $\bar{D}=\underset{v_j\in F}{\sum}  b_{kj}y^{r_{jk}+r}D_{v_j}$.

 If $r_1,\dots,r_h$ \ are the no semisimple  roots associated to  $v_1$ \  and \  $\bar{r}_{ij}=r_i+r_{1j}$ ($r_{11}=0$), then $\bar{r}_{1j},\dots,\bar{r}_{hj}$  are all no semisimple root associated to $v_j$.

We  take the following basis in $V_F$:

$y^{\bar{r}_{11}}D_{v_1},\dots,y^{\bar{r}_{1l}}D_{v_l},y^{\bar{r}_{21}}D_{v_1},\dots,
y^{\bar{r}_{2l}}D_{v_l},\dots,y^{\bar{r}_{h1}}D_{v_1},\dots,y^{\bar{r}_{hl}}D_{v_l}$

The matrix associated to  $\tau_C$ respect to this basis is:

\bigskip

\hskip 4cm $\begin{pmatrix}B'&0&\dots&0\\ 0&B'&\dots&0\\ \vdots&\vdots&\dots&\vdots\\ 0&0&\dots&B'\end{pmatrix}$ \ \  being $B'=(C^{-1})^t$

\medskip

b) Let  $r$ be a no semisimple root associated to $v_k\notin F$ \  and \   $v_F(r)=v_1(r)+\dots+v_l(r)=n$. We can identified $r$ with the element \  $y_1^{v_1(r)}\cdot y_2^{v_2(r)}\cdot{\dots}\cdot y_l^{v_l(r)}\in S^n E_F$. For each $s_1,\dots s_l\in{\mathbb N}$ such that $s_1+\dots+s_l=n$, $r'=r+\overset{l}{\underset{i=2}{\sum}} (s_i-v_i(r))r_{1i}$ \  is a root associated to $v_k$ \  such that $v_i(r')=s_i$.

If we denote \  $V_k^n$ \  the $k$-vectorial space generated by \  $\{y^rD_{v_k}\}_{v_F(r)=n}$, we can identify \  $V_k^n$ \  with \  $S^nE_F$ \  and we have that \   $V_k={\underset{n\geq 0}{\oplus}} V_k^n$ \  where \  $V_k^n=0$ if there is not   no semisimple root $r$ such that $v_F(r)=n$.

 Let \ $C=(c_{ij})\in Gl_F$ \  and \  $D= y^r D_{v_k}$ ($v_k\notin F$).

 $\tau_C\circ\tau_{tD}\circ\tau_{C^{-1}}(y_i)=y_i$  \  for each $y_i\neq y_k$.

$\tau_C\circ\tau_{tD}\circ\tau_{C^{-1}}(y_k)=\tau_C\circ\tau_{tD}(y_k)=\tau_C(y_k+ty_k\cdot y^r)=y_k+t\tau_C(y_k\cdot y^r)$.

Let \  $n_i=v_i(r)$ . $y_k\cdot y^r=y_1^{n_1}\cdot{\dots}\cdot y_l^{n_l}\cdot y$ \  where $y$ is a product of variables different from $y_k$ and  $y_i\in F$.

$\tau_C\circ\tau_{tD}\circ\tau_{C^{-1}}(y_k)=y_k+t(c_{11}y_1+\dots+c_{1l}y_l)^{n_1}\cdot{\dots}\cdot(c_{l1}y_1+\dots+c_{ll}y_l)^{n_l}\cdot y=$

The polynomial \  $(c_{11}y_1+\dots+c_{1l}y_l)^{n_1}\cdot(c_{21}y_1+\dots+c_{2l}y_l)^{n_2}
\cdot{\dots}\cdot(c_{l1}y_1+\dots+c_{ll}y_l)^{n_l}$ \ is a homogeneous polynomial in the variables $y_1,\dots,y_l$ \  of  degree $n=n_1+\dots+n_l$. That is:  the monomials occurring in the sum are the type  \  $y^s=y_1^{s_1}\cdot{\dots}\cdot y_l^{s_l}$ \  where  $s_1+\dots+s_l=n$. Therefore \  $y^s\cdot y =y^s\cdot\dfrac{y_k\cdot y^r}{y_1^{n_1}\cdot{\dots}\cdot y_l^{n_l}}\overset{y^{r_i}=\frac{y_i}{y_1}}{=}y_k\cdot y^{r+(s_2-n_2)r_{12}+\dots+(s_l-n_l)r_{1l}}$

   Therefore \  $\tau_C\circ\tau_{tD}\circ\tau_{C^{-1}}=\tau_{\overline{D}}$ \  where $\overline{D}\in V_k^n$ \  and  we can conclude.

c) With the above notations, if $v_k\nless v_j$ for all $v_j\in F$, then $v_1(r)=\dots=v_l(r)=0$  \ for all no semisimple root associated to $v_k$  and \  $\tau_C\circ\tau_{tD}\circ\tau_{C^{-1}}(y_k)=\tau_{tD}(y_k)$.

\end{pf}

\end{subsection}
\begin{thm} a)  $Aut_g^0A_C=R_u\rtimes (Gl_{F_1}\times\dots\times Gl_{F_k})$ \  where $R_u=E_1\rtimes\dots\rtimes E_l$.

b) The subgroups $L_i=E_1\rtimes\dots\rtimes E_i$ are normal subgroups of $Aut_g^0A_C$.

c) The radical of $Aut_g^0A_C=R_u\rtimes T_{{\mathbb Z}^k}$
\end{thm}

\begin{pf} a) and b) are a consequence of  the theorems 6.6, 6.14 and 6.17.

c) This is because the radical of $Gl_F$ is its center $T_{\mathbb Z}$.
\end{pf}

\section{The group of automorphisms of X}
\begin{subsection}{\bf $Aut^0_kX$}

\

\medskip

We know that \  $Aut^0_kX=Aut^0_gA_C \, / T_{N}=(R_u\rtimes G_R) \, / \, T_{N}$.

\begin{prop} The subgroup $T_{N}\subset Aut^0_gA_C $ \  acts trivially  on $R_u$.
\end{prop}

\begin{pf} We know that \  $0\to M\overset{v}{\to} {\mathbb Z}^r\to N\to 0$ \  is an exact sequence. Therefore \  $T_{N}=\{\lambda \in T_{{\mathbb Z}^r} \, / \lambda^{v(\alpha)}=1 \text{ for all } \alpha\in M\}$.

A root associated to $v_i$  has  $N$-degree $0$. Therefore if \  $D=t\cdot y^rD_{v_i}$, then \   $r=v(\alpha)$ \  for some $\alpha\in M$. Let  \  $\lambda=(\lambda_1,\lambda_2,\dots,\lambda_r)\in T_{N}\subset T_{{\mathbb Z}^r}$.

 $h_\lambda\circ \tau_{tD}\circ h_\lambda^{-1}(y^s)=h_\lambda\tau_{tD}(\lambda^{-s} y^s)=h_\lambda(\lambda^{-s}y^s(1+t\cdot y^r)^{v_i(s)})=$

 $\lambda^{-s}\lambda^sy^s(1+t\lambda^{v(\alpha)}y^r)^{v(s)}=y^s(1+t\cdot y^r)^{v(s)}=
 \tau_{tD}(y^s)$.

\end{pf}

\begin{thm} a) $Aut^0_k X\simeq R_u\rtimes \overline{G}_R$.

b) $L_i$ are normal subgroups of  \ $Aut^0_kX$.

c) The radical of $Aut^0_k X\simeq R_u\rtimes T_{\bar{M}}$ \ where $\bar{M}=M/M'$ and $M'$ is the submodule generated by the semisimple roots.

\end{thm}
\begin{pf} a) and b) are a consequence of theorem 6.18.

c) The radical of \ $\overline{G}_R$ \ is \ $T_{{\mathbb Z}^k}/T_N$. We can conclude because the following sequences are exact:
 $|F_i|=l_i$, $p_i:{\mathbb Z}^{l_i}\to {\mathbb Z}$ is $p_i(n_1,\dots,n_{l_i})=n_1+\dots+n_{l_i}$ \  and

 $Ker\, p_i=<(1,0,\dots,0,-1,0,\dots,0)>$

$$\begin{matrix} &0 & &0 &  & &\\&\downarrow & &\downarrow &  & &\\ 0\to &M' & \overset{v}{\to}&Ker\,p_1\times\dots\times Ker\,p_{k} &\to & 0 & &\\ &\downarrow & &\downarrow &  &\downarrow &\\ 0\to &M&\overset{v}{\to}&{\mathbb Z}^r={\mathbb Z}^{l_1}\times\dots\times {\mathbb Z}^{l_k}&\to& N&\to 0\\
                    & \downarrow& &\hskip 0.6cm\downarrow \hskip 0.5cm p \hskip 0.5cm \downarrow& &\| & \\
0\to &\bar{M}&{\to}&{\mathbb Z}^k={\mathbb Z}\times\dots\times {\mathbb Z}&\to& N&\to 0\\
&\downarrow & &\downarrow & &  & \\
&0 & & 0& &  &\end{matrix}$$
\end{pf}

\end{subsection}
\begin{subsection}{\bf $Aut_kX/Aut^0_kX$}

\

\medskip

Let $X$ and $X'$ be two toric varieties with torus $T_M$ and $T_{M'}$ respectively.
\begin{dfn}A toric morphism from $X$ to $X'$ is a $k$-scheme morphism

$f:X\to X'$ \  and a group morphism \  $\tau: T_M\to T_{M'}$ \  such that \  $f(\lambda\cdot x)=\tau(\lambda)\cdot f(x)$ \  for all $x\in X$ and $\lambda\in T_M$. That is to say \  $f\circ h_\lambda=h_{\tau(\lambda)}\circ f$ \  for each $\lambda\in T_M$.

A toric isomorphism between $X$ and $X'$ is a toric morphism such that the corresponding morphisms of schemes and groups are isomorphism.

We will denote $Aut_{to}X$  to the group of  toric automorphisms of $X$.
\end{dfn}

\begin{rem}Every group morphism   \  $f: T_M\to T_{M'}$ \  is induced by a morphism of ${\mathbb Z}$-module \ $\tau: M'\to M$. The morphism on the ring $k[M']$ is:  $f(x^{m'})=x^{\tau(m')}$ for all $m'\in M'$.

\end{rem}

An example of toric automorphism of $X$ is any element of $T_M$. Indeed $T_M$ is a normal subgroup of $Aut_{to}X$.

Other example is an automorphism of $M$ leaving invariant of the fan:
  Let \  $\tau:M\to M$ \ be  an automorphism and $\tau^*:M^*\otimes_{{\mathbb Z}}\mathbb Q\to M^*\otimes_{{\mathbb Z}}\mathbb Q$ its transpose morphism. Suppose  \ that  for each polyhedral cone   $\sigma\in \Delta$, $\tau^*(\sigma)\in \Delta$. This says that $\tau^{-1}(S_\sigma)=S_{\tau^*(\sigma)}$ and so  $\tau$ induces a toric  automorphism  \  $f_\tau:X\to X$ \  given  by $f_\tau(x^\alpha)=x^{\tau(\alpha)}$ \  for each $\alpha\in M$.

  In particular \  $\tau^*(v_i)=v_{p(i)}$ \  for some permutation $p$ of $\{1,\dots,r\}$ and so there is a finite number of these  automorphisms of this type.

 We denote the set of the   automorphisms of this type by $Aut_\Delta M$.

 \begin{thm} The  natural morphism \  $Aut_\Delta M \to Aut_{to}X/T_M$ \  is an isomorphism.
 \end{thm}

  \begin{pf} Let $ \tau:X\to X$ be a toric isomorphism with automorphism associated $\phi:T_M\to T_M$. If $H$ is a $T_M$-invariant hypersuperface, so is $\tau(H)$. Therefore $\tau(H_i)=H_{p i}$ \  for some permutation $p$ of $\{1,2,\dots,r\}$ and  \  $\tau(U_M)=U_M$.

    Let\  $e\in Spec\,k[M]=U_M$ \  the identity element of the group \  and \  $\tau(e)=\lambda\in T_M$. $\tau'=h_\lambda^{-1}\circ \tau$ \  is a toric isomorphism such that $\tau'(e)=e$. One can check that $\tau'$ restricted to $U_M$ is $\phi$. We also denote by $\phi$ the corresponding automorphism of $M$.
 If $\sigma$ is a polyhedral cone of the fan generated by $v_{i_1},\dots,v_{i_l}$, then  $\tau(Spec\,k[S_\sigma])=Spec\,k[S_{\phi^*(\sigma)}]$ is a $T_M$-invariant affine open of $X$. Therefore $\phi^*(\sigma)$ is a polyhedral cone of the fan. (Theorem 4.2 in [Od]).

\end{pf}

\begin{prop} With the previous notations, if $f:X\to X'$ is an isomorphism of $k$-schemes, there exists   $\tau\in Aut^0_k X'$ and  a toric isomorphism $\bar{f}:X\to X'$ such that $\bar{f}=\tau\circ f$.
\end{prop}
\begin{pf}  We have a group isomorphism: $\phi: Aut_k^0X\to Aut_k^0X'$ \  given by \  $\phi(g)=f\circ g\circ f^{-1}$. As
 $\phi(T_M)$ and $T_{M'}$ \  are maximal torus of \ $Aut_k^0X'$, they are conjugated ([Hu]). That is to say: there exists  \  $\tau\in Aut_k^0X'$ \  such that \  $\tau\circ\phi(T_M)\circ \tau^{-1}=T_{M'}$. Therefore \  $\tau\circ f:X\to X'$ \  is a toric isomorphism with isomorphism associated $h_\tau\circ \phi$ \ where $h_\tau$ is to conjugate by $\tau$.
\end{pf}
\begin{cor}Two  complete toric varieties are isomorphic as algebraic varieties if  and only if they are isomorphic as toric varieties.
\end{cor}

For each equivalence class $F_i$, we denote $S_i=Biy F_i\subset Gl_{F_i}$ the symmetric group of order $|F_i|$ (permutations of the variables belonging to $F_i$). Taking quotient by $T_N$, we have that $S_i$ is a subgroup of $Aut_k^0X$. We show that each $S_i$ is a subgroup of $Aut_\Delta M$. Let $F_i=\{v_1,\dots,v_l\}$ and $p\in S_i$ a permutation of $\{1,2,\dots,l\}$. As automorphism of $A_C$ is such that $\tau_p(y_i)=y_{pi}$ ($pj=j$ for $v_j\notin F_i$).

This produces an automorphism $\tau$ of $k(X)$. $x^\alpha$ is identified with $y^{v(\alpha)}$. If $v(\alpha)=(n_1,n_2,\dots,n_r)$, $\tau(x^\alpha)=\tau_p(y^{v(\alpha)})=\tau_p(y_1^{n_1}\cdot y_2^{n_2}\cdot{\dots}\cdot y_r^{n_r})=y_{p1}^{n_1}\cdot y_{p2}^{n_2}\cdot{\dots}\cdot y_{pr}^{n_r}$. Let $\beta_i$ be the semisimple root associated $(v_i,v_{pi})$ and  $\beta_j=0$ if $v_j\notin F_i$. We have that $y_{pi}=y_iy^{v(\beta_i)}$. Therefore $\tau_p(y^{v(\alpha)})=y^{v(\alpha)}\cdot  y^{n_1v(\beta_1)}\cdot{\dots}\cdot y^{n_rv(\beta_r)}=y^{v(\alpha+n_1\beta_1+\dots+n_r\beta_r)}$. We conclude that $\tau(x^\alpha)=x^{\alpha+n_1\beta_1+\dots+n_r\beta_r}\in k[M]$ and so the automorphism produced in $X$ is the induced by the morphism $M\to M$ given by \  $\alpha\to \alpha+v_1(\alpha)\beta_1+\dots+v_r(\alpha)\beta_r$.

\begin{thm}  $Aut_kX/Aut^0_kX\simeq Aut_\Delta M/(S_1\times\dots\times S_k)$.
\end{thm}

\begin{pf} The natural morphism $Aut_\Delta M\to Aut_kX/Aut^0_kX$ \  is an epimorphism applying theorem 7.5 and proposition 7.6.

 We only have to prove that the kernel of this morphism  is $S_1\times\dots\times S_k$.

Let $\tau: M\to M$ an automorphism such that $\tau^*(v_i)=v_{p(i)}$ for some permutation $p$ of $\{1,\dots,r\}$. As toric automorphism  $\tau:X\to X$, this says that $\tau(H_i)=H_{pi}$. If $\tau\in Aut^0_kX$, by the corollary 2.3,  $H_i\equiv \tau(H_i)=H_{pi}$ and so $p\in S_1\times\dots\times S_k$. $\tau\circ\tau_{p}^{-1}$ is an automorphism in $Aut_\Delta M$ which transforms $v_i$ in $v_i$ and so it is the identity.
\end{pf}

\end{subsection}

\centerline{REFERENCES}

\bigskip

 [Co] Conrad,B. \textit{Reductive Group Schemes} in Autour des Sch\'{e}mas en Groupes Vol. I, Panoramas et
Synth\`{e}ses Numero 42-43, Soci\'{e}t\'{e} Math. de France, 2014.

 [C] Cox, D.A. \textit{The Homogeneous Coordinate Ring of a Toric Variety}. J. Algebraic Geometry 4 (1995), 17-50.

 [C]  Cox, D.A. \textit{Erratum to "The Homogeneos Coordinate Ring of a Toric Variety"}, J.
Algebraic Geometry (2014).

 [D] Demazure M. \textit{Sous-groupes alg\'{e}briques de rang maximun du groupe de Cremona}
Annales scientifiques de l'\'E.N.S. 4e s\'{e}rie, tome 3, n4 (1970), p. 507-588.

 [DG] Demazure M., Gabriel P. \textit{Groupes Alg\'{e}briques}.  Masson, Paris, 1970.

 [F]  Fulton  W. \textit{Introduction to Toric Varieties}. Annals of Mathematics Studies. Number 131. Princeton University Press.

 [H] Hartshorne R. \textit{Algebraic Geometry}. Graduated Texts in Mathematics 52. Springer-Verlag.

 [Hu] Humphreys J. \textit{Linear Algebraic Groups}. Graduated Texts in Mathematics 21. Springer-Verlag.

 [MO] Matsumura H., Oorts F. \textit{Representability of Group Functors and Automorphism of Algebraic Schemes}. Inventiones math 4, 1-25 (1967).

 [Od] Oda T.\textit{Lectures on Torus Embeddings and Applications}. TATA,1978.

 [Sa] Sancho de Salas, C. 2001. Grupos algebraicos y teor\'ia de invariantes, volume 16 of
Aportaciones Matem\'{a}ticas: Textos. Sociedad Matem\'{a}tica Mexicana,M\'{e}xico.

 [SGA3] Demazure M., Grothendieck  A. \textit{ Sch\'{e}mas en groupes I, II, III} Lecture Notes
in Math 151, 152, 153, Springer-Verlag, New York (1970).

\end{document}